\documentclass[twoside,11pt,final]{entics} 
\usepackage{enticsmacro}
\usepackage{graphicx}
\newif\ifignore 
\ignorefalse
\newcommand{\auxproof}[1]{
\ifignore\mbox{}\newline
\textbf{PROOF:} \dotfill\newline
{\it #1}\mbox{}\newline
\textbf{ENDPROOF}\dotfill
\fi}

\newcommand{\ignore}[1]{}

\usepackage{amsmath}
\usepackage{amssymb}
\usepackage{amstext}
\usepackage{mathtools}
\usepackage{mathrsfs}
\usepackage{stmaryrd}
\usepackage{xspace}
\usepackage{enumitem}
\usepackage{url}
\usepackage{nicefrac}
\usepackage{fancybox}
\usepackage{hyperref}
\usepackage{listings}
\usepackage{wrapfig}
\usepackage{pdfrender}
\usepackage{xypic}
\usepackage[all,cmtip]{xy}
\xyoption{2cell}
\UseTwocells
\xyoption{tips}

\usepackage{tikzit}
\usetikzlibrary{circuits.ee.IEC}
\newif\ifexternalizetikz

\newenvironment{myproof}{\begin{trivlist} \item[\hskip \labelsep%
{\bf Proof.}]}{\end{trivlist}}

\newcommand{\QEDbox}{\square}
\newcommand{\QED}{\hspace*{\fill}$\QEDbox$}

\DeclareSymbolFont{T1op}{T1}{cmr}{m}{n}
\SetSymbolFont{T1op}{bold}{T1}{cmr}{bx}{n}
\DeclareMathSymbol{\mathguilsinglleft}{\mathopen}{T1op}{'016}
\DeclareMathSymbol{\mathguilsinglright}{\mathclose}{T1op}{'017}

\newcommand{\flrn}{\ensuremath{\mathsl{Flrn}}}

\newcommand{\concat}{\ensuremath{\mathbin{+{\kern-.5ex}+}}}
\newcommand{\coefm}[1]{\ensuremath{\fatten[0.6pt]{(}{\kern1pt}#1{\kern1pt}\fatten[0.6pt]{)}}}

\newcommand{\flip}{\ensuremath{\mathsl{flip}}}
\newcommand{\binomialname}{\ensuremath{\mathsl{bn}}}
\newcommand{\bivbinomialname}{\ensuremath{\mathsl{bvbn}}}
\newcommand{\binomial}[1][]{\ensuremath{\binomialname[#1]}}
\newcommand{\bivbinomial}[1][]{\ensuremath{\bivbinomialname[#1]}}

\newcommand{\mulnom}{\ensuremath{\mathsl{mn}}}
\newcommand{\multinomial}[1][]{\ensuremath{\mulnom[#1]}}

\newcommand{\poissonname}{\mathsl{pois}}
\newcommand{\poisson}[1][]{\ensuremath{\poissonname[#1]}}

\newcommand{\bibinom}[2]{\left({\kern-.7ex}\binom{#1}{#2}{\kern-.7ex}\right)}
\newcommand{\smallbibinom}[2]{\big({\kern-.4ex}\binom{#1}{#2}{\kern-.4ex}\big)}
\newcommand{\Beta}{\ensuremath{\mathsl{Beta}}}
\newcommand{\Dirichlet}{\ensuremath{\mathsl{Dir}}}

\newcommand{\Dir}{\ensuremath{\mathrm{Dir}}}

\newcommand{\setin}[3]{\{#1\in#2\;|\;#3\}}

\newcommand{\bigsetin}[3]{\big\{#1\in#2\;\big|\;#3\big\}}
\newcommand{\Bigset}[2]{\Big\{{\kern.2em}#1\;\Big|\;#2{\kern.2em}\Big\}}
\newcommand{\Bigsetin}[3]{\Big\{{\kern.2em}#1\in#2\;\Big|\;#3{\kern.2em}\Big\}}

\newcommand{\supp}{\mathrm{supp}}
\newcommand{\mean}{\mathrm{mean}}
\newcommand{\var}{\mathrm{var}}
\newcommand{\cov}{\mathrm{cov}}

\newcommand{\plusje}{\ensuremath{{\kern-2pt}+{\kern-2pt}}}
\newcommand{\minnetje}{\ensuremath{{\kern-1.5pt}-{\kern-1.5pt}}}

\newcommand{\marginalheads}{\ensuremath{\mathsl{heads}}}

\newcommand{\DKL}{\ensuremath{\mathsl{D}_{\mathsl{KL}}}}

\newcommand{\ket}[1]{\ensuremath{|{\kern.1em}#1{\kern.1em}\rangle}}
\newcommand{\bigket}[1]{\ensuremath{\big|{\kern.1em}#1{\kern.1em}\big\rangle}}
\newcommand{\ketstrut}{\vrule height 8.5pt depth 4.5pt width 0pt}
\newcommand{\Bigket}[1]{\ensuremath{\left|\ketstrut{\kern.1em}\right.{\kern-.2em}#1{\kern-.2em}\left.\ketstrut{\kern0em}\right>}}

\newcommand{\zero}{\ensuremath{\mathbf{0}}}
\newcommand{\andthen}{\mathrel{\&}}

\newcommand{\pull}{\mathrel{\mathchoice%
   {\scalebox{-0.5}[1]{$\gg=$}}
   {\scalebox{-0.5}[1]{$\gg{\kern-1.5ex}=$}}
   {\scalebox{-0.5}[1]{${\kern.5ex}\scriptstyle\gg{\kern-0.2ex}={\kern.5ex}$}}
   {\scalebox{-0.5}[1]{$\scriptscriptstyle\gg=$}}}}
\newcommand{\push}{\mathrel{\mathchoice%
   {\scalebox{-0.5}[1]{$=\ll$}}
   {\scalebox{-0.5}[1]{$={\kern-1.5ex}\ll$}}
   {\scalebox{-0.5}[1]{${\kern.5ex}\scriptstyle={\kern-0.2ex}\ll{\kern.5ex}$}}
   {\scalebox{-0.5}[1]{$\scriptscriptstyle=\ll$}}}}

\newcommand{\pushing}[2]{{#1}_{*}(#2)}

\newcommand{\distributionsymbol}{\mathcal{D}}

\newcommand{\multisetsymbol}{\mathcal{M}}

\newcommand{\Dst}{\distributionsymbol}

\newcommand{\Mlt}{\multisetsymbol}

\newcommand{\UF}{\ensuremath{\mathcal{U}{\kern-.75ex}\mathcal{F}}}

\newcommand{\Kl}{\mathcal{K}{\kern-.4ex}\ell}
\newcommand{\EM}{\mathcal{E}{\kern-.4ex}\mathcal{M}}

\newcommand{\NNO}{\mathbb{N}}
\newcommand{\pNNO}{\mathbb{N}_{>0}}

\newcommand{\R}{\mathbb{R}}
\newcommand{\nnR}{\R_{\geq 0}}

\newcommand{\finset}[1]{\ensuremath{\boldsymbol{#1}}}

\newcommand{\Ef}{\ensuremath{\mathcal{E}{\kern-.5ex}f}}
\newcommand{\intd}{{\kern.2em}\mathrm{d}{\kern.03em}}

\newcommand{\OF}{\ensuremath{\mathcal{O}{\kern-.1em}\mathcal{F}}}
\newcommand{\Closed}{\ensuremath{\mathcal{C}{\kern-.05em}\ell}}

\newcommand{\congrightarrow}{\mathrel{\smash{\stackrel{
           \raisebox{.5ex}{$\scriptstyle\cong$}}{
           \raisebox{0ex}[0ex][0ex]{$\rightarrow$}}}}}

\newcommand{\Prob}{\footnotesize \mathrm{Pr}}

\usetikzlibrary{shapes,shapes.geometric}
\usetikzlibrary{shapes.gates.ee,shapes.gates.ee.IEC}
\usetikzlibrary{calc}
\usetikzlibrary{positioning}
\usetikzlibrary{cd}
\usetikzlibrary{decorations.pathmorphing}

\ifexternalizetikz
\usepackage{tikzcdx}
\usetikzlibrary{external}
\makeatletter
\newcommand{\tikzextname}[1]{%
\tikzset{external/figure name={\tikzexternal@realjob-#1-}}}
\makeatother
\fi

\newsavebox\sbpto
\savebox\sbpto{\begin{tikzpicture}[baseline=-2.5pt]
            \filldraw[fill=white,draw=white] circle (1.4pt);
            \filldraw[fill=white,draw=black,line width=0.2pt]circle(2pt);
                \end{tikzpicture}}
\newcommand\chanto{\mathrel{\ooalign{$\to$\cr
      \hfil\!$\usebox\sbpto$\hfil\cr}}}            
          
\newsavebox\sbground
\savebox\sbground{\begin{tikzpicture}[circuit ee IEC,yscale=1,xscale=1]
                \draw (0,-2ex) to (0,0) node[ground,rotate=90,xshift=.65ex] {};
                \end{tikzpicture}}

\newsavebox\sbunif
\savebox\sbunif{\begin{tikzpicture}[circuit ee IEC,yscale=1,xscale=1]
                \draw (0,0) to (0,2ex) node[ground,rotate=270,xshift=2.5ex] {};
                \end{tikzpicture}}

\tikzset{dot/.style =
  {inner sep=0mm,minimum width=1mm,minimum height=1mm,
    draw,shape=circle}}
\tikzset{minicopy/.style = {dot,fill,text depth=-0.2mm}}
\newsavebox\sbcopier
\savebox\sbcopier{%
  \begin{tikzpicture}[baseline=0pt]
    \node[minicopy,scale=.7] (a) at (0,3.6pt) {};
    \draw (a) -- +(-90:.30);
    \draw (a) -- +(45:.35);
    \draw (a) -- +(135:.35);
  \end{tikzpicture}}

\newsavebox\sbcocopier
\savebox\sbcocopier{%
  \begin{tikzpicture}[baseline=0pt]
    \node[minicopy,scale=.7] (a) at (0,3.6pt) {};
    \draw (a) -- +(90:.30);
    \draw (a) -- +(-45:.35);
    \draw (a) -- +(-135:.35);
  \end{tikzpicture}}

\newsavebox\sbcup
\savebox\sbcup{%
  \begin{tikzpicture}[baseline=0pt]
    \node (a) at (-10pt,18pt) {};
    \node (b) at (10pt,18pt) {};
    \draw (a) to[out=-90,in=-90,looseness=2] (b);
  \end{tikzpicture}}

\newsavebox\sbcap
\savebox\sbcap{%
  \begin{tikzpicture}[baseline=0pt]
    \node (a) at (-10pt,-6pt) {};
    \node (b) at (10pt,-6pt) {};
    \draw (a) to[out=90,in=90,looseness=2] (b);
  \end{tikzpicture}}

\ifexternalizetikz\tikzexternalenable\fi

\newcommand*{\fatten}[1][.4pt]{%
  \textpdfrender{
    TextRenderingMode=FillStroke,
    LineWidth={\dimexpr(#1)\relax},
  }%
}

\ifdefined\mathsl\else
  \DeclareMathAlphabet{\mathsl}{\encodingdefault}{\rmdefault}{\mddefault}{\sldefault}
  \SetMathAlphabet{\mathsl}{bold}{\encodingdefault}{\rmdefault}{\bfdefault}{\sldefault}
\fi

\makeatletter
\newcommand{\mathoverlap}[2]{\mathpalette\mathoverlap@{{#1}{#2}}}
\newcommand{\mathoverlap@}[2]{\mathoverlap@@{#1}#2}
\newcommand{\mathoverlap@@}[3]{\ooalign{$\m@th#1#2$\crcr\hidewidth$\m@th#1#3$\hidewidth}}
\makeatother

\newcommand\inline[1]{{\lstinline!#1!}}

\def\conf{MFPS 2025} 	
\volume{5}			
\def\volu{5}			
\def\papno{10}			
\def\mydoi{17042}
\def\lastname{Jacobs} 
\def\runtitle{A fresh look at Bivariate Bionomial Distributions} 
\def\copynames{B. Jacobs}    

\begin{document}

\begin{frontmatter}

\title{A Fresh Look at Bivariate Binomial Distributions} 						
  \author{Bart Jacobs\thanksref{myemail}}	
   \address{iHub \\ Radboud University\\				
    Nijmegen, The Netherlands}  						
   \thanks[myemail]{Email: \href{mailto:bart@cs.ru.nl} {\texttt{\normalshape
        bart@cs.ru.nl}}} 

\begin{abstract} 
Binomial distributions capture the probabilities of `heads' outcomes
when a (biased) coin is tossed multiple times. The coin may be
identified with a distribution on the two-element set $\{0,1\}$, where
the $1$ outcome corresponds to `head'. One can also toss two separate
coins, with different biases, in parallel and record the
outcomes. This paper investigates a slightly different `bivariate'
binomial distribution, where the two coins are dependent (also called:
entangled, or entwined): the two-coin is a distribution on the product
$\{0,1\} \times \{0,1\}$. This bivariate binomial exists in the
literature, with complicated formulations. Here we use the language of
category theory to give a new succint formulation.  This paper
investigates, also in categorically inspired form, basic properties of
these bivariate distributions, including their mean, variance and
covariance, and their behaviour under convolution and under updating,
in Laplace's rule of succession. Furthermore, it is shown how
Expectation Maximisation works for these bivariate binomials, so that
mixtures of bivariate binomials can be recognised in data. This paper
concentrates on the bivariate case, but the binomial distributions may
be generalised to the multivariate case, with multiple dimensions, in
a straightforward manner.
\end{abstract}

\end{frontmatter}

\section{Introduction}\label{IntroSec}

A very basic physical model of a random phenomenon is the tossing (or
`flipping' or `throwing') of a coin. The possible outcomes are `head'
and `tail', or `$1$' and `$0$' in this paper. The coin may have a bias
$r\in [0,1]$, which we incorporate in a `flip' distribution on the set
of outcomes $\{0,1\}$. Using `ket' notation $\ket{-}$ it is written
as:
\begin{equation}
  \label{FlipEqn}
  \begin{array}{rcl}
  \flip(r)
  & = &
  r\bigket{1} + (1-r)\bigket{0}.
  \end{array}
\end{equation}

\noindent This says that the probability of outcome $1$ is $r$ and the
probability of outcome $0$ is $1-r$.

\begin{figure}
\begin{center}
  \includegraphics[scale=0.3]{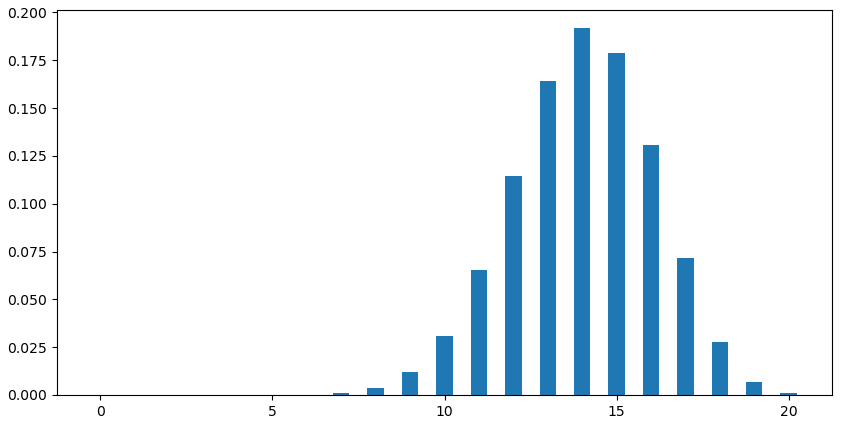}
  \hspace*{4em}
  \includegraphics[scale=0.5]{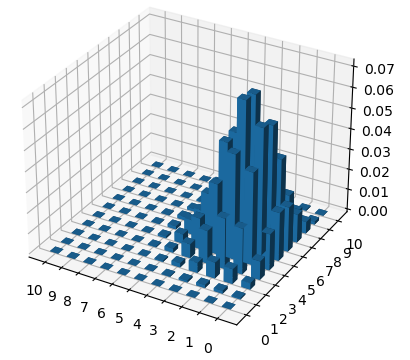}
\end{center}
\caption{Plots of an ordinary binomial distribution (20 tosses with
  bias $\frac{7}{10}$) on the left, and a bivariate binomial
  distribution on the right, with 10 tosses of a two-coin
  $\frac{3}{8}\bigket{0, 0} + \frac{5}{12}\bigket{0, 1} +
  \frac{1}{12}\bigket{1, 0} + \frac{1}{8}\bigket{1, 1}$.}
\label{BinomialsFig}
\end{figure}

Repeating such coin flipping $K\in\NNO$ times, still with a bias $r\in
[0,1]$, yields the binomial distribution on the set $\{0,1,\ldots,K\}$
where a number $n$ in the set gets the probability that $n$ out of
these $K$ flips are head. Again using ket notation this
binomial distribution is written as:
\begin{equation}
  \label{BinomialEqn}
  \begin{array}{rcl}
  \binomial[K](r)
  & = &
  \displaystyle\sum_{0\leq n\leq K} \, \binom{K}{n} \cdot r^{n}\cdot (1-r)^{K-n}
  \,\bigket{n}.
  \end{array}
\end{equation}

\noindent The probabilities add up to one by the Binomial Theorem.  On
the left in Figure~\ref{BinomialsFig} one sees the plot of such a
binomial distribution with $K=20$ tosses of a coin with bias $r =
\frac{7}{10}$. Via the bias this typical (discrete) bell shape can be
shifted to the left or right.

On the right in Figure~\ref{BinomialsFig} we see a plot of a bivariate
binomial distribution. It has a bell shape in two dimensions. It does
not arise from a single coin distribution $\flip(r)$ on $\{0,1\}$ but
from what we call a two-coin distribution on
$\{0,1\}\times\{0,1\}$. Such a two-coin may arise as a tensor product
$\flip(r)\otimes\flip(s)$ of two separate coins, but in general it is
more than such a product. It may incorporate the important feature of
`joint' distributions on product spaces, namely that the whole
distribution does not coincide with the product of its marginals.  In
that case the joint distribution is called dependent, entangled, or
entwined.

\begin{wrapfigure}[9]{r}{0pt}
\includegraphics[scale=0.2]{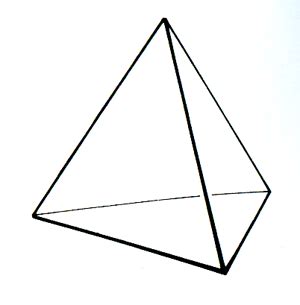}
\end{wrapfigure}
A possible physical model of such a two-coin distribution on
$\{0,1\}\times\{0,1\} = \{(0,0), (0,1), (1,0), (1,1)\}$ is given by a
tetrahedron, as on the side, also known as a triangular pyramid. Its
four faces can be labeled $00$, $01$, $10$ and $11$. Such a
tetrahedron can be tossed, where the outcome is the label on the face
on which the tetrahedron lands. One could imagine that the weight
distribution inside the tetrahedron is such that the outcomes are
determined by a general two-coin distribution:
\[ r_{1}\bigket{0,0} + r_{2}\bigket{0,1} + r_{3}\bigket{1,0} + r_{4}\bigket{1,1}
\qquad\mbox{where}\qquad
r_{1} + r_{2} + r_{3} + r_{4} = 1. \]

This paper formulates this bivariate binomial distribution, for such a
two-coin, using some basic category theory, namely the functoriality
of sending a set $X$ to the set $\Dst(X)$ of distributions on $X$. For
a number $K\in\NNO$ of tosses, it produces a distribution on the
product set $\{0,1,\ldots,K\} \times \{0,1,\ldots,K\}$ --- as on the
right in Figure~\ref{BinomialsFig} where $K=10$.  This bivariate
binomial is introduced as a (functorial) pushforward of a multinomial
distribution along a `marginal heads' function, see
Section~\ref{BivBinBinomSec} for the mathematical
details. Intuitively, one throws the tetrahedron $K$-many times and
one separately counts the number of $1$'s in the first $a$ and in the
second $b$ position in the label $ab$ written on the faces that the
tetrahedron lands on (in $K$-many tosses).  We focus on a
$2$-dimensional extension of the binomial distribution. One can also
extend to $N$-dimensional form, via a `hedron' object with $2^N$ sides
and labels from $\{0,1\}^{N}$. It gives a distribution on
$\{0,1,\ldots,K\}^{N}$ when one counts the $1$'s in each of the $N$
positions $1,\ldots,N$ of the (down-facing) labels, after multiple
tosses.

The (ordinary) binomial distribution~\eqref{BinomialEqn} can be seen
as discrete version of the widely used Gaussian (normal) distribution.
Gaussians often occur in $N$-dimensional form. The $N$-dimensional
discrete version, however, is less familiar. In Section~\ref{EMSec} we
show how to recognise (mixtures of) bivariate binomials in (multisets
of) data items. This technique may have practical value and could
increase the popularity of bivariate binomial distributions.

This paper starts by describing the relevant background information,
on multisets and distributions, tensor products and convolutions of
distributions, channels and their inversions, and on binomial and
multinomial distributions. Section~\ref{BivBinBinomSec} then
introduces bivariate binomial distributions, via functoriality, and
establishes some basic properties with respect to marginalisation and
also closure under convolution. Subsequently, Section~\ref{MeanSec}
determines the mean, variance and covariance of bivariate
binomials. It turns out that they can be expressed, respectively, as
$K$ times the mean, variance and covariance of the underlying two-coin
distribution, where $K$ is the number of tosses.
Section~\ref{LaplaceSec} briefly discusses Laplace's rule of
succession, which is in generally used to describe the expected
outcome after an update. Here it is formulated in general terms, as a
mean of a dagger (channel), and more specifically for bivariate
distributions, with Dirichlet and Poisson distributions as
priors. Section~\ref{EMSec} gives an illustration of how bivariate
binomials can be recognised in suitable data multisets.

\section{Preliminaries}\label{PrelimSec}

\subsection{Multisets, distributions and random variables}\label{MltDstSubsec}

A \emph{multiset} is like a subset, except that its elements may occur
multiple times. Alternatively, a multiset is like a list in which the
order of the elements does not matter. For instance, an urn filled
with three red, two green and five blue balls forms a multiset. We
write it using ket notation as $3\ket{R} + 2\ket{G} + 5\ket{B}$, where
the ket brackets $\ket{-}$ serve to separate the elements $R,G,B$ for
the three colours, from their multiplicities $3,2,5$. Also a draw of
multiple balls from such an urn forms a multiset. In general, multiple
data items, as used for learning (see Section~\ref{EMSec}) are
naturally described as a multiset.

For an arbitrary set $X$ we shall write $\Mlt(X)$ for the set of
(finite) multisets with elements from $X$. A multiset
$\varphi\in\Mlt(X)$ may be described in two equivalent ways, namely
via kets as a formal sum $\varphi = n_{1}\ket{x_1} + \cdots +
n_{K}\ket{x_{K}} = \sum_{i} n_{i}\ket{x_i}$, with elements $x_{i} \in
X$ and associated multiplicities $n_{i}\in\NNO$. Alternatively, we may
describe a multiset $\varphi \in\Mlt(X)$ as a function $\varphi\colon
X \rightarrow \NNO$ with finite support, where $\supp(\varphi) =
\setin{x}{X}{\varphi(x) \neq 0}$. We can combine these notations and
write $\varphi = \sum_{x} \varphi(x)\ket{x}$, where the summands with
multiplicity zero are omitted. 

The size $\|\varphi\|\in\NNO$ of a multiset $\varphi\in\Mlt(X)$ is the
total number of its elements, including multiplicities. Thus,
$\|\varphi\| = \sum_{x} \varphi(x)$. The abovementioned urn $\upsilon
= \ket{R} + 2\ket{G} + 5\ket{B}$ has size $\|\upsilon\| = 10$.  There
is one multiset in $\Mlt(X)$ of size zero, namely the empty multiset,
written as $\zero\in\Mlt(X)$. As a function, it satisfies $\zero(x) =
0$, for each $x\in X$, giving $x$ multiplicity~$0$.  We write
$\Mlt[K](X) \subseteq \Mlt(X)$ for the subset of multisets of
size $K\in\NNO$.

The mapping $X \mapsto \Mlt(X)$ forms a monad on the category of
sets. We only use functoriality: a function $f\colon X \rightarrow
Y$ yields another function $\Mlt(f) \colon \Mlt(X) \rightarrow \Mlt(Y)$
via the following two equivalent definitions, in ket form on the left,
and in function form on the right.
\begin{equation}
  \begin{array}{rclcrcl}
  \label{MltFunctorEqn}
  \Mlt(f)\Big(\sum_{i} n_{i}\bigket{x_i}\Big)
  & = &
  \sum_{i} n_{i}\bigket{f(x_{i})}
  & \hspace*{4em} &
  \Mlt(f)(\varphi)(y)
  & = &
  \displaystyle\sum_{x \in f^{-1}(y)} \, \varphi(x).
  \end{array}
\end{equation}

\noindent We shall especially use functoriality for a projection
function $\pi_{1} \colon X \times Y \rightarrow X$. Then
$\Mlt(\pi_{1}) \colon \Mlt(X\times Y) \rightarrow \Mlt(X)$ performs
marginalisation.

One more observation is that $\Mlt(X)$ is the free commutative monoid
on the set $X$. The monoid structure is given by pointwise addition of
multisets: $\big(\varphi+\psi\big)(x) = \varphi(x) + \psi(x)$, with
the everywhere zero multset $\zero\in\Mlt(X)$ as neutral element.
Functions $\Mlt(f)$ preserve this monoid structure.

We turn to distributions: a (finite discrete probability) distribution
$r_{1}\ket{x_{1}} + \cdots + r_{K}\ket{x_K}$ is written as a multiset,
except that the multiplicities $r_{i}$ are now probabilities from the
unit interval $[0,1]$ that add up to one: $\sum_{i} r_{i} = 1$. We
shall write $\Dst(X)$ for the set of distributions with elements from
$X$. A distribution $\omega\in\Dst(X)$ may be written in ket form, as
we just did, but also as a function $\omega \colon X \rightarrow
[0,1]$ with finite support $\supp(\omega) = \setin{x}{X}{\omega(x)
  \neq 0}$ and with $\sum_{x} \omega(x) = 1$.  When $X$ is a finite
set, we say that $\omega\in\Dst(X)$ has full support when
$\supp(\omega) = X$, that is, when $\omega(x) > 0$ for each $x\in
X$. The operation $\Dst$ is also a monad on the category of
sets. Functoriality works as for $\Mlt$, in~\eqref{MltFunctorEqn}, in
particular with projections $\Dst(\pi_{i})$ for marginalisation.

Each non-empty multiset $\varphi\in\Mlt(X)$ can be turned into a
distribution $\flrn(\varphi)\in\Dst(X)$, where $\flrn$ stands for
frequentist learning. This involves learning by counting:
\begin{equation}
\label{FlrnEqn}
\begin{array}{rcl}
\flrn(\varphi) & \coloneqq & \displaystyle\sum_{x\in\supp(\varphi)} \,
\frac{\varphi(x)}{\|\varphi\|}\,\bigket{x}.
\end{array}
\end{equation}

\noindent For instance, for the above urn $\upsilon = 3\ket{R} +
2\ket{G} + 5\ket{B}$ we get $\flrn(\upsilon) = \frac{3}{10}\ket{R} +
\frac{1}{5}\ket{G} + \frac{1}{2}\ket{B}$. It captures the
probabilities of drawing a ball of a particular colour from the urn.

A random variable is a pair $(\omega, p)$ where $\omega\in\Dst(X)$ is
a distribution and $p\colon X \rightarrow \R$ an `observable'. The
validity (or expected value) is:
\begin{equation}
  \label{ValidityEqn}
  \begin{array}{rcl}
  \omega\models p
  & \,\coloneqq\, &
  \displaystyle\sum_{x\in X} \, \omega(x) \cdot p(x).
  \end{array}
\end{equation}

\subsection{Products and convolutions of distributions}\label{ProdConvSubsec}

For two distributions $\omega\in\Dst(X)$ and $\rho\in\Dst(Y)$ one can
form the (parallel) product distribution $\omega\otimes\rho \in
\Dst(X\times Y)$, via $\big(\omega\otimes\rho\big)(x,y) =
\omega(x)\cdot\rho(y)$. Marginalisation of such a product yields
the components:
\[ \begin{array}{rclcrcl}
  \Dst(\pi_{1})\big(\omega\otimes\rho\big)
  & = &
  \omega
  & \qquad\mbox{and}\qquad &
  \Dst(\pi_{2})\big(\omega\otimes\rho\big)
  & = &
  \rho.
\end{array} \]

\noindent In general, for a `joint' distribution $\tau \in
\Dst(X\times Y)$, on a product space, one does \emph{not} have $\tau =
\Dst(\pi_{1})(\tau) \otimes \Dst(\pi_{2})(\tau)$. If this equation
does hold, $\tau$ is called independent, non-entangled, or
non-entwined.

We write $\finset{2} = \{0,1\}$. A distribution $\tau \in
\Dst(\finset{2}\times\finset{2})$, say $\tau = r_{1}\ket{0,0} +
r_{2}\ket{0,1} + r_{3}\ket{1,0} + r_{4}\ket{1,1}$ with $r_{1} + r_{2}
+ r_{3} + r_{4} = 1$, is non-entwined if and only if $r_{1}\cdot r_{4}
= r_{2}\cdot r_{3}$.

Tensors of distributions are used for a special sum of distributions,
called convolution. Let $M = (M,+,0)$ be a commutative monoid and let
$\omega,\rho\in\Dst(M)$ be two distributions on the underlying set
$M$. Then we can form their convolution $\omega+\rho \in\Dst(M)$ as,
via functoriality of $\Dst$, applied to the addition $+ \colon M\times
M \rightarrow M$, as in:
\begin{equation}
  \label{ConvolutionEqn}
  \begin{array}{rcccl}
    \omega + \rho
    & \coloneqq &
    \Dst(+)\Big(\omega\otimes\rho\Big)
    & = &
    \displaystyle\sum_{x,y\in M} \, \omega(x)\cdot\rho(y)\,\bigket{x+y}.
  \end{array}
\end{equation}

\noindent This turns $\Dst(M)$ into a commutative monoid, with
$1\ket{0} \in\Dst(M)$ as unit, for $0\in M$. It is not hard to see
that if $f\colon M \rightarrow M'$ is a homomorphism of monoids, then
so is $\Dst(f) \colon \Dst(M) \rightarrow \Dst(M')$.

In the sequel we shall encounter distributions on $\{0,1\}$ or more
generally on $\{0,\ldots,K\}$. By considering these sets as subsets of
$\NNO$, with its additive commutative monoid structure, we can apply
the convolution definition~\eqref{ConvolutionEqn} to such
distributions. Commutative monoids are closed under products $\times$,
so that $\NNO^{2} = \NNO\times\NNO$ is also a commutative monoid, via
component-wise sums and zeros. Also, as we have seen, the set
$\Mlt(X)$ of multisets over a set $X$ is a commutative monoid, so that
convolution can be used for distributions on multisets. This applies
to multinomial distributions, see Subsection~\ref{BinomMulnomSubsec}.

\subsection{Channels and their inversions / daggers}\label{ChanSubsec}

A channel is a function of the form $c\colon X \rightarrow
\Dst(Y)$. It maps an element $x\in X$ to a distribution $c(x)$ on
$Y$. As such it occurs often as conditional probability, written as
$\Prob(y|x)$. Channels are well-behaved probabilistic computations
that can be composed sequentially and in parallel. They form the
morphisms in the symmetric monoidal Kleisli category $\Kl(\Dst)$ of
the distribution monad $\Dst$. In fact we simplify the notation and
use a special arrow $\chanto$ for channels. Thus we simply write
$c\colon X \chanto Y$ for the more cumbersome $c\colon X \rightarrow
\Dst(Y)$.

A basic operation is pushforward along a channel. Given a
channel $c\colon X \chanto Y$ and a distribution $\omega\in\Dst(X)$ on
its domain $X$, we can form a new distribution $\pushing{c}{\omega}
\in \Dst(Y)$ on its codomain $Y$, namely:
\begin{equation}
\label{PushEqn}
\begin{array}{rcccl}
\pushing{c}{\omega}
& \coloneqq &
\displaystyle\sum_{x\in X} \omega(x)\cdot c(x)
& = &
\displaystyle\sum_{y\in Y} \, 
   \left(\sum_{x\in X} \omega(x)\cdot c(x)(y)\right)\bigket{y}.
\end{array}
\end{equation}




Under suitable circumstances a channel can be turned around, via an
operation called Bayesian inversion, see
\textit{e.g.}~\cite{ClercDDG17,ChoJ19,Fritz20}. This involves turning
a conditional probability $\Prob(y|x)$ into $\Prob(x|y)$.  We present
it in channel form: so let $c\colon X \chanto Y$ be a channel from $X$
to $Y$, where the set $Y$ is finite. We will turn it into a channel $Y
\chanto X$, in the opposite direction. This inversion requires a
`prior' distribution $\omega\in\Dst(X)$ such that the pushforward
$\pushing{c}{\omega} \in \Dst(Y)$ has full support. Then we can define
a new channel, written as $c^{\dag}_{\omega} \colon Y \chanto X$ via:
\begin{equation}
\label{DaggerEqn}
\begin{array}{rcl}
c^{\dag}_{\omega}(y)
& \coloneqq &
\displaystyle\sum_{x\in X} \, \frac{\omega(x)\cdot c(x)(y)}
   {\pushing{c}{\omega}(y)} \, \ket{x}.
\end{array}
\end{equation}

\noindent The full support requirement is needed to prevent
division-by-zero errors.

\subsection{Binomial and multinomial distributions}\label{BinomMulnomSubsec}

We have seen the binomial distribution $\binomial[K](r) \in
\Dst\big(\{0,\ldots,K\}\big)$ in Equation~\eqref{BinomialEqn}, for
$K$-many tosses of a coin with bias $r\in [0,1]$. A single coin toss
is given by the flip distribution $\flip(r)$ in~\eqref{FlipEqn}.

\begin{lemma}
\label{BinomialLem}
\begin{enumerate}
\item \label{BinomialLemConv} The binomial distribution is an iterated
  convolution of the flip distribution:
    \[ \begin{array}{rcccl}
      \binomial[K](r)
      & = &
      K\cdot\flip(r)
      & = &
      \underbrace{\flip(r) + \cdots + \flip(r)}_{K\text{ times}}.
    \end{array} \]

    \noindent Hence, binomial distributions are closed under
    convolution: $\binomial[K](r) + \binomial[L](r) =
    \binomial[K+L](r)$.

\item \label{BinomialLemMeanVar} The mean and variance of the binomial
  distribution are given by:
\[ \begin{array}{rclcrcl}
\mean\Big(\binomial[K](r)\Big)
& = &
K\cdot r
& \qquad\mbox{and}\qquad &
\var\Big(\binomial[K](r)\Big)
& = &
K\cdot (K-1)\cdot r.
\end{array} \eqno{\QEDbox} \]
\end{enumerate}
\end{lemma}

We turn to the multinomial distribution. We can view a distribution
$\omega\in\Dst(X)$ as an abstract urn, where $X$ is the set of colours
and $\omega(x)\in [0,1]$ gives the probability of drawing a ball of
colour $x\in X$. We are interested in drawing multiple balls, in
`multinomial mode', with replacement. A draw of $K$-many balls from
the urn $\omega$ can be identified with a $K$-sized multiset
$\varphi\in\Mlt[K](X)$. The multinomial distribution
$\multinomial[K](\omega) \in \Dst\big(\Mlt[K](X)\big)$ assigns
probabilites to such draws of size $K$. It is defined as:
\begin{equation}
  \label{MultinomialEqn}
  \begin{array}{rclcrcl}
\multinomial[K](\omega)
& \coloneqq &
\displaystyle\sum_{\varphi\in\Mlt[K](X)} \, \coefm{\varphi} \cdot 
\prod_{x\in \supp(\omega)} \omega(x)^{\varphi(x)}\, \bigket{\varphi}
& \quad\mbox{where}\quad &
\coefm{\varphi}
& \coloneqq &
\displaystyle\frac{\|\varphi\|!}{\prod_{x} \varphi(x)!}.
  \end{array}
\end{equation}

\noindent The precise form of this multinomial distribution does not
matter so much at this stage; interested readers may
consult~\cite{Jacobs21b,Jacobs22a}. We make the following facts
explicit. The first two items are standard. For details about the
third item, see \textit{e.g.}~\cite[Prop.~7.2]{Jacobs25a}
or~\cite{Jacobs21g}.

\begin{lemma}
  \label{MultinomialLem}
Consider the multinomial distribution~\eqref{MultinomialEqn}
for $\omega\in\Dst(X)$ and $K\in\NNO$.
  \begin{enumerate}
  \item \label{MultinomialLemBin} The binomial and multinomial
    distribution can be identified when $X = \finset{2} = \{0,1\}$,
    via the isomorphisms $\flip \colon [0,1] \congrightarrow
    \Dst(\finset{2})$ and $\{0,\ldots,K\} \congrightarrow
    \Mlt[K](\finset{2})$, via $n\mapsto n\ket{1} + (K-n)\ket{0}$, in:
    \[ \xymatrix@R-1pc{
      [0,1]\ar[rr]^-{\binomial[K]}\ar@{=}[d]_{\rotatebox{90}{$\sim$}}
      & & \Dst\big(\{0,\ldots,K\}\big)\ar@{=}[d]_{\rotatebox{90}{$\sim$}}
      \\
      \Dst(\finset{2})\ar[rr]^-{\multinomial[K]} & &
      \Dst\big(\Mlt[K](\finset{2})\big)
    } \]

  \item \label{MultinomialLemConv} Multinomials are closed under
    convolution: $\multinomial[K](\omega) + \multinomial[L](\omega) =
    \multinomial[K+L](\omega)$. 

  \item \label{MultinomialLemElts} For elements $y,z\in X$ with
    $y\neq z$,
\[ \begin{array}[b]{rcl}
\displaystyle\sum_{\varphi\in\Mlt[K](X)}\, 
   \multinomial[K](\omega)(\varphi) \cdot \varphi(y)
& = &
K \cdot \omega(y)
\\[+1em]
\displaystyle\sum_{\varphi\in\Mlt[K](X)}\, 
   \multinomial[K](\omega)(\varphi) \cdot \varphi(y) \cdot \varphi(z)
& = &
K \cdot (K - 1) \cdot \omega(y) \cdot \omega(z).
\\[+1em]
\displaystyle\sum_{\varphi\in\Mlt[K](X)}\, 
   \multinomial[K](\omega)(\varphi) \cdot \varphi(y)^{2}
& = &
K \cdot (K - 1) \cdot \omega(y)^{2} \,+\, K\cdot\omega(y).
\end{array} \eqno{\QEDbox} \]
  \end{enumerate}
\end{lemma}

\section{Introducing bivariate distributions}\label{BivBinBinomSec}

This section introduces the main distribution of this paper, namely
the bivariate binomial distribution. It also describes some of its
most basic properties. Our definition exploits functoriality of
distributions $\Dst$, via a marginal heads function that is defined
first. This definition is compared to explicit formulations in the
literature in Remark~\ref{FormulationRem} below.

\begin{definition}
\label{MarginalHeadsDef}
Fix a number $K\in\NNO$ and consider the following `marginal heads'
function:
\[ \xymatrix{
\Mlt[K]\big(\finset{2}\times \finset{2}\big)\ar[r]^-{\marginalheads} &
   \{0,1,\ldots,K\} \times \{0,1,\ldots,K\},
} \]

\noindent defined as pair of marginals evaluated at $1 \in
\finset{2}$, in:
\[ \begin{array}{rcccl}
  \marginalheads(\varphi)
  & \coloneqq &
  \big(\,\varphi(1,0) + \varphi(1,1), \, \varphi(0,1) + \varphi(1,1)\,\big)
  & = &
  \big(\, \Mlt(\pi_{1})(\varphi)(1), \,
  \Mlt(\pi_{2})(\varphi)(1) \,\big).
\end{array} \]

\noindent The last line makes clear that the multiplicities of $1$,
that is of `head', in the marginal mulitsets are counted by this
function $\marginalheads$.

More generally, for a `dimension' $N\geq 1$ one can define:
\begin{equation}
  \label{NaryMarginalHeadsEqn}
  \vcenter{\xymatrix@R-2.2pc{
  \Mlt[K]\big(\finset{2}^{N}\big)\ar[rr]^-{\marginalheads}
  & & \{0,\ldots,K\}^{N}
  \\
  \varphi\ar@{|->}[rr] & & \big(\, \Mlt(\pi_{1})(\varphi)(1),
  \ldots, \Mlt(\pi_{N})(\varphi)(1) \,\big)
  }}
\end{equation}
\end{definition}

We shall need the following combinatorial result about the marginal
heads function.

\begin{lemma}
\label{MarginalHeadsCoefficientLem}
  For $K\geq 0$ and $N\geq 1$ and numbers $\vec{n} = n_{1}, \ldots,
  n_{N} \in \{0,\ldots,K\}$ one has:
  \[ \begin{array}{rcl}
    \displaystyle\sum_{\varphi\in\Mlt[K](\finset{2}^{N}), \,
      \marginalheads(\varphi) = \vec{n}} \,
    \coefm{\varphi}
    & \,=\, &
    \displaystyle\prod_{1\leq i \leq N} \, \binom{K}{n_{i}}.
  \end{array} \]
\end{lemma}

\begin{myproof}
  We do the proof for $N=2$ and use Vandermonde's formula, see
  \textit{e.g.}~\cite{GrahamKP94}: for numbers $B,G\in\NNO$ and $K\leq
  B+G$,
\begin{equation}
\label{VandermondeBinaryEqn}
\begin{array}{rcl}
\displaystyle\binom{B+G}{K}
& = &
\displaystyle\sum_{b\leq B, \, g\leq G, \, b+g = K}\, 
   \binom{B}{b}\cdot\binom{G}{g}.
\end{array}
\end{equation}
  
\noindent Now, for a multiset
$\varphi\in\Mlt[K]\big(\finset{2}\times\finset{2}\big)$ we use new
variables $i,j$ below, where $i = \varphi(1,1)$ and $j =
\varphi(0,1)$.
\[ \begin{array}[b]{rcl}
  \lefteqn{\sum_{\varphi\in\Mlt[K](\finset{2}\times\finset{2}), \,
      \marginalheads(\varphi) = (n_{1}, n_{2})} \, \coefm{\varphi}}
  \\
  & = &
  \displaystyle\sum_{\varphi\in\Mlt[K](\finset{2}\times\finset{2}), \,
    \varphi(1,0) + \varphi(1,1) = n_{1},\,\varphi(0,1) + \varphi(1,1) = n_{2}}
    \frac{K!}{\varphi(0,0)!\cdot\varphi(1,0)!\cdot\varphi(0,1)!
      \cdot\varphi(1,1)!}
    \\[+1em]
    & = &
  \displaystyle\sum_{i\leq n_{1}, \, j\leq K-n_{1}, \, i+j=n_{2}} \,
  \frac{K!}{(K- n_{1} - j)! \cdot (n_{1}- i)!
    \cdot j! \cdot i!}
  \\[+0.8em]
  & = &
  \displaystyle\sum_{i\leq n_{1}, \, j\leq K-n_{1}, \, i+j=n_{2}} \,
  \frac{K!}{n_{1}!\cdot (K- n_{1})!} \cdot
  \frac{n_{1}!}{i!\cdot (n_{1} - i)!} \cdot
  \frac{(K- n_{1})!}{j! \cdot (K- n_{1} - j)!}
  \\[+1.2em]
  & = &
  \displaystyle\binom{K}{n_1} \cdot
  \sum_{i\leq n_{1}, \, j\leq K-n_{1}, \, i+j=n_{2}} \,
  \binom{n_1}{i} \cdot \binom{K- i}{j}
  \hspace*{\arraycolsep}\;\smash{\stackrel{\eqref{VandermondeBinaryEqn}}{=}}\;\hspace*{\arraycolsep}
  \displaystyle\binom{K}{n_{1}} \cdot \binom{K}{n_{2}}.    
\end{array} \eqno{\QEDbox} \]
\end{myproof}

\begin{definition}
\label{BivBinDef}
For a `two-coin' distribution $\gamma\in\Dst\big(\finset{2} \times
\finset{2}\big)$ and a `toss number' $K$ the \emph{bivariate} binomial
distribution is defined as pushforward along the marginal heads
function, of a multinomial distribution:
\begin{equation}
\label{BivBinEqn}
\begin{array}{rcl}
\bivbinomial[K](\gamma)
& \coloneqq &
\Dst\big(\marginalheads\big)\Big(\multinomial[K](\gamma)\Big)
   \,\in\, \Dst\Big(\{0,1,\ldots,K\} \times \{0,1,\ldots,K\}\Big)
\\
& = &
\displaystyle \sum_{0\leq n_{1}, n_{2} \leq K} \;\;
   \sum_{\varphi\in\Mlt[K](\finset{2}\times\finset{2}), \,
   \marginalheads(\varphi) = (n_{1}, n_{2})} \,
   \multinomial[K](\gamma)(\varphi)\,\bigket{n_{1}, n_{2}}.
\end{array}
\end{equation}

\noindent The $N$-dimensional version of this binomial distribution is
obtained by first forming the multinomial distribution
$\multinomial[K](\gamma) \in
\Dst\big(\Mlt[K]\big(\finset{2}^{N}\big)\big)$, for an $N$-coin
$\gamma\in\Dst\big(\finset{2}^{N}\big)$, and then pushing along the
$N$-dimensional marginal heads function~\eqref{NaryMarginalHeadsEqn}.
\end{definition}

\begin{example}
  \label{BivBinEx}
Consider the two-coin distribution $\gamma = \frac{3}{8}\bigket{0, 0}
+ \frac{5}{12}\bigket{0, 1} + \frac{1}{12}\bigket{1, 0} +
\frac{1}{8}\bigket{1, 1}$. It is not hard to see that it is entwined,
\textit{i.e.}~that it is different from the product $\otimes$ of its
marginals. We use $K=2$ and first consider the multinomial
distribution on $\Mlt[2]\big(\finset{2}\times\finset{2}\big)$. It is
written via `kets of kets' in:
\[ \begin{array}{rcl}
  \multinomial[2](\gamma)
  & = &
  \frac{9}{64}\Bigket{2\ket{0, 0}} +
  \frac{5}{16}\Bigket{1\ket{0, 0} + 1\ket{0, 1}} +
  \frac{25}{144}\Bigket{2\ket{0, 1}} +
  \frac{1}{16}\Bigket{1\ket{0, 0} + 1\ket{1, 0}}
  \\
  & & \quad +\,
  \frac{5}{72}\Bigket{1\ket{0, 1} + 1\ket{1, 0}} +
  \frac{1}{144}\Bigket{2\ket{1, 0}} +
  \frac{3}{32}\Bigket{1\ket{0, 0} + 1\ket{1, 1}}
  \\
  & & \quad +\,
  \frac{5}{48}\Bigket{1\ket{0, 1} + 1\ket{1, 1}} +
  \frac{1}{48}\Bigket{1\ket{1, 0} + 1\ket{1, 1}} +
  \frac{1}{64}\Bigket{2\ket{1, 1}}.
\end{array} \]

\noindent We now obtain the associated bivariate binomial
distribution by applying the marginals heads function $\marginalheads$
to the multisets inside the above big kets $\Bigket{-}$. By
counting the total number of $1$'s on the left and on the right
in these multisets one obtains:
\[ \begin{array}{rcl}
  \bivbinomial[2](\gamma)
  \hspace*{\arraycolsep}=\hspace*{\arraycolsep}
  \Dst\big(\marginalheads\big)\Big(\multinomial[2](\gamma)\Big)
  & = &
  \frac{9}{64}\bigket{0, 0} +
  \frac{5}{16}\bigket{0, 1} +
  \frac{25}{144}\bigket{0, 2} +
  \frac{1}{16}\bigket{1, 0} +
  \frac{47}{288}\bigket{1, 1}
  \\[+0.2em]
  & & \quad + \,
  \frac{5}{48}\bigket{1, 2} +
  \frac{1}{144}\bigket{2, 0} +
  \frac{1}{48}\bigket{2, 1} +
  \frac{1}{64}\bigket{2, 2}.
\end{array} \]

\noindent The bivariate binomial for draw size $K=10$, of the
same $\gamma$, is plotted on the right in
Figure~\ref{BinomialsFig}. 
\end{example}

\begin{remark}
\label{FormulationRem}
\begin{enumerate}
\item \label{FormulationRemOther} Our formulation of the bivariate
  binomial in Definition~\ref{BivBinDef} uses the functoriality of
  $\Dst$.  Common definitions in the literature do not exploit this
  functoriality and involve more complicated explicit formulations.
  For instance, in~\cite[Eqn.~(1)]{BayramogluE10} one finds the
  following expression, translated into ket form:
\begin{equation}
\label{BivBinAltEqn}
\begin{array}{l}
\displaystyle \sum_{0\leq k,l\leq K} \; \sum_{\max(0, k+l-K) \leq i \leq \min(k,l)} \;
   \frac{K!}{i!\cdot (k-i)! \cdot (l-i)! \cdot (K-k-l+i)!} 
\\[-0.4em]
\hspace*{15em} \cdot
   \gamma(0,0)^{i} \cdot \gamma(0,0)^{k-i} \cdot 
   \gamma(0,0)^{l-i} \cdot \gamma(0,0)^{K-k-l+i} \,\Bigket{k,l}
\end{array}
\end{equation}

\noindent One recognises elements of the multinomial
distribution~\eqref{MultinomialEqn} in this formulation, but it is
hard to follow what is going on. In fact, this
formulation~\eqref{BivBinAltEqn} is subtly different, from the one
that we use in~\eqref{BivBinEqn}, since it does not count \emph{heads}
/ \emph{ones} but it counts \emph{tails} / \emph{zeros}, in each
coordinate separately.

Multi-dimensional (multivariate) versions of the
formula~\eqref{BivBinAltEqn} become even harder, see
\textit{e.g.}~\cite{BayramogluE10,Teugels90}. Here, the
multi-dimensional form is obtained in a straightforward manner: the
definition~\eqref{BivBinEqn} still works for an $N$-coin
$\gamma\in\Dst\big(\finset{2}^{N}\big)$ with the marginal heads
function~\eqref{NaryMarginalHeadsEqn} used in $N$-ary form.

\ignore{

# https://dm.ieu.edu.tr/Ism/ozge.pfd

TwoTwo = range_sp(2) @ range_sp(2)

def bivbin(n, twocoin):
    out_sp = range_sp(n+1) ** 2
    return DState(lambda k, l:
                  sum([ math.factorial(n) / (math.factorial(i) *
                                             math.factorial(k-i) *
                                             math.factorial(l-i) *                                                         math.factorial(n-k-l+i)) *
                        twocoin(0,0) ** i *
                        twocoin(0,1) ** (k-i) *
                        twocoin(1,0) ** (l-i) *
                        twocoin(1,1) ** (n-k-l+i)
                      for i in range(max(0, k+l-n), min(k,l)+1) ]),
                  out_sp)

w = random_distribution(TwoTwo)
n = 3
print("")
b = bivbin(4, w)
print( b.state_size() )
print( b )
print("")
print( BivariateBinomial_functorial(4)(w) )

}

\item \label{FormulationRemOwn} We can also give a more explicit
  formulation of our functorial definition~\eqref{BivBinEqn} of the
  bivariate binomial, in the style of~\eqref{BivBinAltEqn}, via the
  following characterisation of the inverse image.
\begin{equation}
\label{MarginalHeadsInverseEqn}
\begin{array}{rcl}
\lefteqn{\bigsetin{\varphi}{\Mlt[K]\big(\finset{2}\times\finset{2}\big)}
   {\marginalheads(\varphi) = (n_{1}, n_{2})}}
\\[+0.2em]
& = &
\begin{cases}
\begin{array}{l}
\big\{(K - n_{2} - i)\bigket{0,0} + (n_{2} - n_{1} + i)\bigket{0,1} + 
   i\bigket{1,0} + (n_{1} - i)\bigket{1,1} 
\\
\qquad\qquad \Big|\;0 \leq i \leq \min(n_{1}, K-n_{2})\big\}
\end{array}  & \mbox{if } n_{1} \leq n_{2}
\\[+1.4em]
\begin{array}{l}
\big\{(K - n_{1} - i)\bigket{0,0} + i\bigket{0,1} + 
   (n_{1} - n_{2} + i)\bigket{1,0} + (n_{2} - i)\bigket{1,1} 
\\
\qquad\qquad \Big|\;0 \leq i \leq \min(n_{2}, K-n_{1})\big\}
\end{array}  & \mbox{if } n_{2} < n_{1}.
\end{cases}
\end{array} 
\end{equation}

\ignore{

def BivariateBinomial(K):
    """Optimised Bivariate Binomial distribution """
    out_sp = range_sp(K+1) @ range_sp(K+1)
    def chan_fun(dist):
        sp = dist.sp
        zero = zero_frac if dist.frac else 0
        mn = Multinomial(K)(dist)
        def state_fun(n1,n2):
            if n2 <= n1:
                return list_addition(
                    [ mn(DState([K-n1-i, i, n1-n2+i, n2-i], sp))
                      for i in range(min(n2, K-n1) + 1) ],
                    initialiser = zero)
            else:
                return list_addition(
                    [ mn(DState([K-n2-i, n2-n1+i, i, n1-i], sp))
                      for i in range(min(n1, K-n2) + 1) ],
                    initialiser = zero)
        return DState(state_fun, out_sp, frac=dist.frac)
    return AChannel(chan_fun, 1)

}

\end{enumerate}
\end{remark}

We first collect some basic, standard properties about the bivariate
binomial distribution, in particular in relation to the ordinary
binomial distribution. The definition~\eqref{BivBinEqn} of the bivariate
binomial distribution via functoriality allows equational reasoning
in proofs. The details are relegated to the appendix.

\begin{lemma}
\label{BivBinMarginalLem}
Let a number $K\in\NNO$ and a two-coin distribution
$\gamma\in\Dst\big(\finset{2} \times \finset{2}\big)$ be given.  For
clarity we use different notation for the various projections:
\begin{equation}
  \label{ProjDiag}
  \vcenter{\xymatrix{
    \finset{2}\times\finset{2}\ar@<+0.5ex>[r]^-{\pi_1}\ar@<-0.5ex>[r]_-{\pi_2}
    & \finset{2}
    & &
    \{0,\ldots, K\}\times\{0,\ldots, K\}
    \ar@<+0.5ex>[r]^-{\Pi_1}\ar@<-0.5ex>[r]_-{\Pi_2}  & \{0,\ldots K\}
  }}
  \end{equation}

\noindent We then write $\gamma_{1} \coloneqq \Dst(\pi_{1})(\gamma)$
and $\gamma_{2} \coloneqq \Dst(\pi_{2})(\gamma)$ in $\Dst(\finset{2})$
for the two marginal coins.  They are determined by the `heads'
probabilities $\gamma_{1}(1), \gamma_{2}(1) \in [0,1]$, via the
isomorphism $\Dst(\finset{2}) \cong [0,1]$.
\begin{enumerate}
\item \label{BivBinMarginalLemProj} The marginals of the
  bivariate binomial distribution are the (ordinary) binomials
  of the marginals of the two-coin distribution:
\[ \begin{array}{rclcrcl}
\Dst\big(\Pi_{1}\big)\Big(\bivbinomial[K](\gamma)\Big)
& = &
\binomial[K]\Big(\gamma_{1}(1)\Big)
& \qquad\qquad &
\Dst\big(\Pi_{2}\big)\Big(\bivbinomial[K](\gamma)\Big)
& = &
\binomial[K]\Big(\gamma_{2}(1)\Big).
\end{array} \]

\item \label{BivBinMarginalLemTensor} When $\gamma$ is
  non-entwined, that is, when $\gamma = \gamma_{1}\otimes\gamma_{2}$,
  then $\bivbinomial[K](\gamma)$ is also non-entwined, with:
  \[ \begin{array}{rcl}
    \bivbinomial[K]\Big(\gamma_{1}\otimes\gamma_{2}\Big)
    & = &
    \binomial[K]\Big(\gamma_{1}(1)\Big) \otimes 
   \binomial[K]\Big(\gamma_{2}(1)\Big).
  \end{array} \]
\end{enumerate}
\end{lemma}

This second result confirms that the notion of bivariate
binomial is a well-defined extension of the ordinary binomial: if the
two-coin distribution $\gamma$ is a tensor product of two coins, the
bivariate binomial is a tensor product of the two (marginal,
ordinary) binomial distributions.

There is one more well-behavedness property that we add, namely
closure under convolution, as we have seen for binomial and
multinomial distributions, in Lemmas~\ref{BinomialLem}
and~\ref{MultinomialLem}. Again, the proofs are in the appendix.

\begin{proposition}
  \label{BivBinConvolutionProp}
  Let $\gamma\in\Dst\big(\finset{2}\times\finset{2}\big)$ be given,
  with numbers $K,L\in\NNO$.  Bivariate binomials are closed under
  convolution, as expressed on the left below.  
    \[ \begin{array}{rclcrcccl}
      \bivbinomial[K](\gamma) \,+\, \bivbinomial[L](\gamma)
      & = &
      \bivbinomial[K+ L](\gamma)
      & \hspace*{5em} &
      \bivbinomial[K](\gamma)
      & = &
      K\cdot\gamma
      & = &
      \gamma + \cdots + \gamma.
    \end{array} \]

\noindent Thus, the bivariate binomial is a $K$-fold convolution of
the two-coin distribution $\gamma$, as on the right. \QED
\end{proposition}

\section{Mean, variance and covariance}\label{MeanSec}

This section describes common statistical characteristics of the
bivariate binomial distribution, namely the mean, variance and
covariance. It turns out that these values are $K$-times the values
for the underlying two-coin distributions. This gives nice
formulas. The proofs require some work.

First, for a distribution $\omega\in\Dst\big(\R^{N}\big)$, for some
number $N\geq 1$, the mean is described as the tuple in $\R^{N}$ given
by the validities of the projections $\pi_{i} \colon \R^{N}
\rightarrow \R$, viewed as observables. This gives the formula:
\begin{equation}
  \label{MeanEqn}
  \begin{array}{rcl}
    \mean(\omega)
    & \coloneqq &
    \big(\, \omega\models\pi_{1}, \, \ldots, \, \omega\models\pi_{N} \,\big).
  \end{array}
\end{equation}

\begin{lemma}
  \label{BivBinMeanLem}
For $K\in\NNO$ and $\gamma\in\Dst\big(\finset{2}\times\finset{2}\big)$,
\[ \begin{array}{rcl}
\mean\Big(\bivbinomial[K](\gamma)\Big)
& = &
K\cdot\mean\big(\gamma\big) \;\in\; \R\times\R.
\end{array} \]
\end{lemma}

\begin{myproof}
An abstract argument uses that the mean ommutes with convolutions, in
combination with
Proposition~\ref{BivBinConvolutionProp}. But
one can also proceed in a more concrete manner. First, the mean of
$\gamma\in\Dst\big(\finset{2}\times\finset{2}\big)$ can be expressed
in terms of $\gamma$'s marginals $\gamma_{1} \coloneqq
\Dst(\pi_{1})(\gamma), \gamma_{2} \coloneqq \Dst(\pi_{2})(\gamma) \in
\Dst(\finset{2})$, namely as:
\[ \begin{array}{rcl}
  \mean\big(\gamma\big)
  \hspace*{\arraycolsep}\smash{\stackrel{\eqref{MeanEqn}}{=}}\hspace*{\arraycolsep}
  \big(\, \omega\models\pi_{1}, \, \omega\models\pi_{2} \,\big)
  & = &
  \displaystyle\left(
  \sum_{b_{1}, b_{2} \in \finset{2}} \, \gamma(b_{1},b_{2})\cdot b_{1}, \,
  \sum_{b_{1}, b_{2} \in \finset{2}} \, \gamma(b_{1},b_{2})\cdot b_{2} \right)
  \\[+1em]
  & = &
  \displaystyle\left(
  \sum_{b_{2} \in \finset{2}} \, \gamma(1,b_{2}), \,
  \sum_{b_{1} \in \finset{2}} \, \gamma(b_{1},1) \right)
 \hspace*{\arraycolsep}=\hspace*{\arraycolsep}
  \Big(\,\gamma_{1}(1), \, \gamma_{2})(1)\,\Big).
\end{array} \eqno{(*)} \]

\noindent Next, using the projection notation from~\eqref{ProjDiag},
\[ \begin{array}[b]{rcl}
  \mean\Big(\bivbinomial[K](\gamma)\Big)
  & \smash{\stackrel{\eqref{MeanEqn}}{=}} &
  \big(\, \bivbinomial[K](\gamma)\models\Pi_{1}, \, 
   \bivbinomial[K](\gamma)\models\Pi_{2} \,\big)
   \\[+0.2em]
  & = &
  \displaystyle\left(
  \sum_{n_{1}, n_{2} \in \NNO} \, \bivbinomial[K](\gamma)(n_{1},n_{2})\cdot n_{1}, \,
  \sum_{n_{1}, n_{2} \in \NNO} \, \bivbinomial[K](\gamma)(n_{1},n_{2})\cdot n_{2} \right)
  \\[+1em]
  & = &
  \displaystyle\left(
  \sum_{n_{1} \in \NNO} \, \Dst(\Pi_{1})
   \Big(\bivbinomial[K](\gamma)\Big)(n_{1})\cdot n_{1}, \,
  \sum_{n_{2} \in \NNO} \, \Dst(\Pi_{2})
   \Big(\bivbinomial[K](\gamma)\Big)(n_{2})\cdot n_{2} \right)
  \\[+1em]
  & = &
  \displaystyle\left(
  \sum_{n_{1} \in \NNO} \, \binomial[K]\big(\gamma_{1}(1)\big)(n_{1})\cdot n_{1}, \,
  \sum_{n_{2} \in \NNO} \, 
  \binomial[K]\big(\gamma_{2}(1)\big)(n_{2})\cdot n_{2} \right)
    \quad\mbox{by Lemma\rlap{~\ref{BivBinMarginalLem}~\eqref{BivBinMarginalLemProj}}}
  \\[+1em]
  & = &
  \Big(\,\mean\big(\binomial[K]\big(\gamma_{1}(1)\big)\big), \,
  \mean\big(\binomial[K]\big(\gamma_{2}(1)\big)\big) \, \Big)
  \\[+0.4em]
  & = &
  \Big(\, K\cdot \gamma_{1}(1), \, K\cdot \gamma_{2}(1) \, \Big)
     \qquad \mbox{by Lemma~\ref{BinomialLem}~\eqref{BinomialLemMeanVar}}
  \\[+0.2em]
  & = &
  K\cdot \Big(\,\gamma_{1}(1), \, \gamma_{2})(1)\,\Big)
  \;\;\smash{\stackrel{(*)}{=}}\;\;
  \mean\big(\gamma\big).
\end{array} \eqno{\QEDbox} \]
\end{myproof}

Analogous results exist for the variance and covariance of
the bivariate binomial. We first recall the definitions,
for a distribution $\omega\in\Dst(X)$ and two observables $p,q\colon
X \rightarrow \R$ on the same set:
\begin{equation}
\label{VarCovarEqn}
\begin{array}{rcl}
\var\big(\omega, p)
& = &
\Big(\omega \models p\andthen p\Big) - \Big(\omega\models p\Big)^{2}
\\[+0.4em]
\cov\big(\omega, p, q)
& = &
\Big(\omega \models p\andthen q\Big) - 
   \Big(\omega\models p\Big)\cdot \Big(\omega\models q\Big).
\end{array}
\end{equation}

\noindent The conjunction $\andthen$ of observables is given by
point-wise multiplication. When the underlying set $X$ is (a subset
of) $\R^{N}$ one typically computes these variances and covariances
for the projections. That will be done below.

\begin{theorem}
\label{BivBinVarCovarThm}
For $K\in\NNO$ and $\gamma\in\Dst\big(\finset{2}\times\finset{2}\big)$
one has, for $i,j\in\{1,2\}$,
\[ \begin{array}{rclcrcl}
\var\Big(\bivbinomial[K](\gamma), \Pi_{i}\Big)
& = &
K\cdot\var(\gamma, \pi_{i})
\qquad\mbox{and}\qquad
\cov\Big(\bivbinomial[K](\gamma), \Pi_{i}, \Pi_{j}\Big)
& = &
K\cdot\cov(\gamma, \pi_{i}, \pi_{j}).
\end{array} \]
\end{theorem}

\begin{myproof}
We do the case for $i=1$. First,
\[ \begin{array}{rcl}
  \var\Big(\gamma, \pi_{1}\Big)
  \hspace*{\arraycolsep}\smash{\stackrel{\eqref{VarCovarEqn}}{=}}\hspace*{\arraycolsep}
  \Big(\gamma \models \pi_{1}\andthen\pi_{1}\Big) \,-\,
  \Big(\gamma \models \pi_{1}\Big)^{2}
  & = &
  \displaystyle\left(\sum_{b_{1},b_{2}\in\finset{2}} \,
  \gamma(b_{1}, b_{2})\cdot b_{1} \cdot b_{1}\right) \,-\,
  \left(\sum_{b_{1},b_{2}\in\finset{2}} \,
  \gamma(b_{1}, b_{2})\cdot b_{1}\right)^{2}
  \\[+1em]
  & = &
  \gamma(1,0) + \gamma(1,1) \,-\,
  \Big(\gamma(1,0) + \gamma(1,1)\Big)^{2}.      
\end{array} \]

\noindent We turn to the bivariate binomial and make use of
Lemma~\ref{MultinomialLem}~\eqref{MultinomialLemElts}.
\[ \begin{array}{rcl}
  \var\Big(\bivbinomial[K](\gamma), \Pi_{i}\Big)
  & = &
  \Big(\bivbinomial[K](\gamma) \models \Pi_{1}\andthen\Pi_{1}\Big) \,-\,
  \Big(\bivbinomial[K](\gamma) \models \Pi_{1}\Big)^{2}
  \\[+0.2em]
  & = &
  \displaystyle\left(\sum_{n_{1},n_{2}\in\NNO} \,
  \bivbinomial[K](\gamma)(n_{1}, n_{2})\cdot n_{1}\cdot n_{1}\right) \,-\,
  \left(\sum_{n_{1},n_{2}\in\NNO} \,
  \bivbinomial[K](\gamma)(n_{1}, n_{2})\cdot n_{1}\right)^{2}
  \\[+1em]
  & = &
  \displaystyle\left(\sum_{\varphi\in\Mlt[K](\finset{2}\times\finset{2})} \,
  \multinomial[K](\gamma)(\varphi)\cdot
  \big(\varphi(1,0)+\varphi(1,1)\big)^{2}\right)
  \\
  & & \qquad\displaystyle -\,
  \left(\sum_{\varphi\in\Mlt[K](\finset{2}\times\finset{2})} \,
  \multinomial[K](\gamma)(\varphi)\cdot
  \big(\varphi(0,1)+\varphi(1,1)\big)\right)^{2}
  \\[+1.2em]
  & = &
  K\cdot (K- 1) \cdot \gamma(1,0)^{2} + K\cdot\gamma(1,0) \,+\,
  2\cdot K\cdot (K- 1) \cdot \gamma(1,0) \cdot \gamma(1,1)
  \\
  & & \qquad + \,
  K\cdot (K- 1) \cdot \gamma(1,1)^{2} + K\cdot\gamma(1,1) -
  \Big(K\cdot\gamma(1,0) + K\cdot\gamma(1,1)\Big)^{2}
  \\[+0.2em]
  & = &
  K\cdot\Big(-\gamma(1,0)^{2} + \gamma(1,0) -
  2\cdot \gamma(1,0) \cdot \gamma(1,1)
  -\gamma(1,1)^{2} + \gamma(1,1)\Big)
  \\[+0.2em]
  & = &
  K\cdot\Big(\gamma(1,0) + \gamma(1,1) -
  \big(\gamma(1,0) + \gamma(1,1)\big)^{2}\Big)
  \\[+0.2em]
  & = &
  K\cdot \var\Big(\gamma, \pi_{1}\Big) \quad\mbox{as shown above.}
\end{array} \]

For covariance, the cases where $i=j$ are covered by the previous
point, since covariance with equal observables is variance. We do the
case $i=1,j=2$, which gives the same outcome as for $i=2,j=1$. We
first look at the covariance of the two-coin distribution $\gamma$.
\[ \begin{array}{rcl}
  \cov\Big(\gamma, \pi_{1}, \pi_{2}\Big)
  & \smash{\stackrel{\eqref{VarCovarEqn}}{=}} &
  \Big(\gamma \models \pi_{1}\andthen\pi_{2}\Big) \,-\,
  \Big(\gamma \models \pi_{1}\Big)\cdot\Big(\gamma \models \pi_{2}\Big)
  \\[+0.2em]
  & = &
  \displaystyle\left(\sum_{b_{1},b_{2}\in\finset{2}} \,
  \gamma(b_{1}, b_{2})\cdot b_{1} \cdot b_{2}\right) \,-\,
  \left(\sum_{b_{1},b_{2}\in\finset{2}} \,
  \gamma(b_{1}, b_{2})\cdot b_{1}\right) \cdot
  \left(\sum_{b_{1},b_{2}\in\finset{2}} \,
  \gamma(b_{1}, b_{2})\cdot b_{2}\right)
  \\[+1em]
  & = &
  \gamma(1,1) \,-\,
  \Big(\gamma(1,0) + \gamma(1,1)\Big) \cdot
  \Big(\gamma(0,1) + \gamma(1,1)\Big).
\end{array} \]

\noindent Next, again by Lemma~\ref{MultinomialLem}~\eqref{MultinomialLemElts}:
\[ \hspace*{-1em}\begin{array}[b]{rcl}
  \lefteqn{\cov\Big(\bivbinomial[K](\gamma), \Pi_{1}, \Pi_{2}\Big)}
  \\[+0.2em]
  & = &
  \Big(\bivbinomial[K](\gamma) \models \Pi_{1}\andthen\Pi_{2}\Big) \,-\,
  \Big(\bivbinomial[K](\gamma) \models \Pi_{1}\Big)\cdot
  \Big(\bivbinomial[K](\gamma) \models \Pi_{2}\Big)
  \\[+0.2em]
  & = &
  \displaystyle\left(\sum_{n_{1},n_{2}\in\NNO} \,
  \bivbinomial[K](\gamma)(n_{1}, n_{2})\cdot n_{1}\cdot n_{2}\right)
  \\
  & & \quad \displaystyle -\,
  \left(\sum_{n_{1},n_{2}\in\NNO} \,
  \bivbinomial[K](\gamma)(n_{1}, n_{2})\cdot n_{1}\right)\cdot
  \rlap{$\displaystyle\left(\sum_{n_{1},n_{2}\in\NNO} \,
    \bivbinomial[K](\gamma)(n_{1}, n_{2})\cdot n_{2}\right)$}
  \\[+1em]
  & = &
  \displaystyle\left(\sum_{\varphi\in\Mlt[K](\finset{2}\times\finset{2})} \,
  \multinomial[K](\gamma)(\varphi)\cdot
  \big(\varphi(1,0)+\varphi(1,1)\big)\cdot
  \big(\varphi(0,1)+\varphi(1,1)\big)\right)
  \\
  & & \quad -\,
  \Big(K\cdot\gamma(1,0) + K\cdot\gamma(1,1)\Big) \cdot
  \Big(K\cdot\gamma(0,1) + K\cdot\gamma(1,1)\Big)
  \\[+0.2em]
  & = &
  K\cdot (K- 1)\cdot \Big(\gamma(1,0)\cdot\gamma(0,1) +
  \gamma(1,0)\cdot\gamma(1,1) + \gamma(1,1)\cdot\gamma(0,1)
  \\
  & & \quad +\,
  \gamma(1,1)^{2}\Big) + K\cdot\gamma(1,1) \,-\,
  K^{2} \cdot \Big(\gamma(1,0) + \gamma(1,1)\Big) \cdot
  \rlap{$\Big(\gamma(0,1) + \gamma(1,1)\Big)$}
  \\[+0.2em]
  & = &
  K \cdot \Big(\gamma(1,1) \,-\,
  \big(\gamma(1,0) + \gamma(1,1)\big) \cdot
  \big(\gamma(0,1) + \gamma(1,1)\big)\Big)
  \\[+0.2em]
  & = &
  \cov\Big(\gamma, \pi_{1}, \pi_{2}\Big), \quad\mbox{see before.}
\end{array} \eqno{\QEDbox} \]
\end{myproof}

We conclude with an observation which will be useful for
Expectation Maximisation in Section~\ref{EMSec}.

%
%

\begin{fact}
\label{RecoverFact} Suppose we have a distribution $\sigma
  \in \Dst\Big(\{0,\ldots,K\} \times \{0,\ldots,K\}\Big)$, for some
  given $K\in\NNO$. We know that $\sigma$ is a bivariate binomial
  distribution, of the form $\sigma = \bivbinomial[K](\gamma)$ but we
  do not know $\gamma\in\Dst\big(\finset{2}\times\finset{2}\big)$. We
  can compute the mean and the covariance of $\sigma$. Via
  Proposition~\ref{BivBinMeanLem} and Theorem~\ref{BivBinVarCovarThm}
  we can then obtain the probabilities of $\gamma$, and thus $\gamma$
  itself, in the following way.
\begin{itemize}
\item $\gamma(1,1) = \frac{1}{K} \cdot \Big( \cov\big(\sigma, \Pi_{1},
  \Pi_{2}\big) - \big(\sigma\models\Pi_{1}\big) \cdot
  \big(\sigma\models\Pi_{2}\big)\Big)$;

\item $\gamma(1,0) = \frac{1}{K} \cdot \big(\sigma\models\Pi_{1}\big) -
  \gamma(1,1)$;

\item $\gamma(0,1) = \frac{1}{K} \cdot \big(\sigma\models\Pi_{2}\big) -
  \gamma(1,1)$;

\item $\gamma(0,0) = 1 - \gamma(1,0) - \gamma(0,1) - \gamma(1,1)$.
\end{itemize}

\noindent Thus, when the draw size $K$ is given, a bivariate
binomial distribution is entirely determined by its mean and covariance,
since its mean is given by the pair of validities $(\,
\bivbinomial[K](\gamma) \models \Pi_{1}, \, \bivbinomial[K](\gamma)
\models \Pi_{2} \,)$, see~\eqref{MeanEqn}.
\end{fact}

\section{Laplace's rule of succession}\label{LaplaceSec}

Laplace's rule of succession provides an answer to what was originally
formulated as the sunrise problem: suppose I have seen the sun rise in
$n$ successive days, what is the probability that it will rise
tomorrow? The answer according to Laplace's rule is
$\frac{n+1}{n+2}$. It will go to one as $n$ goes to infinity.

The binomial distribution plays a role in a systematic answer
formulation of the rule of succession. Hence the question arises: is
there a role for the bivariate binomial too? In order to answer
this question we first provide a general, modern interpretation of the
rule of succession. It is generally understood as providing the expected 
outcome after an update, see \textit{e.g.}~\cite{Ross18}
or~\cite{BlitzsteinH19}. We reformulate this as: the rule of
succession calculates the \emph{mean of a dagger}.

The relevant calculations involve continuous distributions as priors,
such as the Beta and Dirichlet distributions. We have to assume that
the reader is reasonably familiar with these distributions and with
the associated results. Therefore, this section only gives a sketch of
the situation at hand, from a structural angle. We start by recalling,
without proof, basic facts about the Beta distribution on $[0,1]$ and
about the Dirichlet distribution on discrete distributions
$\Dst(X)$. These continuous distributions are typically used for the
bias parameter $r\in [0,1]$ of a flip distribution $\flip(r)$ and for
the distribution parameter $\omega\in\Dst(X)$ of a multinomial
distribution $\multinomial[K](\omega)$. The results below are
standard. The formulation of the conjugate prior property in terms of
daggers is in line with~\cite{Jacobs20a}.

\begin{lemma}
\label{BetaLem}
We write $\Beta(\alpha,\beta)$ for the Beta distribution on the unit
interval $[0,1]$, where $\alpha,\beta\in\pNNO$. 
\begin{enumerate}
\item $\mean(\Beta(\alpha,\beta)\big) = \frac{\alpha}{\alpha+\beta}$.

\item The Beta distribution on $[0,1]$ is conjugate prior to the
  binomial channel $\binomial[K] \colon [0,1] \chanto \{0,\ldots,K\}$.
  The dagger channel $\binomial[K]^{\dag}_{\Beta(\alpha,\beta)} \colon
  \{0,\ldots,K\} \chanto [0,1]$, with $\Beta(\alpha,\beta)$ as prior
  distribution, gives a parameter update of the form:
\[ \begin{array}{rcl}
\binomial[K]^{\dag}_{\Beta(\alpha,\beta)}(n)
& = &
\Beta\big(\alpha+n, \, \beta+K-n\big).
\end{array} \]

\item Laplace's rule of succession gives in this case as mean of the
dagger:
\[ \begin{array}{rcccl}
\mean\Big(\binomial[K]^{\dag}_{\Beta(\alpha,\beta)}(n)\Big)
& = &
\mean\Big(\Beta\big(\alpha+n, \, \beta+K-n\big)\Big)
& = &
\displaystyle\frac{\alpha+n}{\alpha+\beta+K}.
\end{array} \eqno{\QEDbox} \]
\end{enumerate}
\end{lemma}

Starting from a uniform prior on $[0,1]$, for $\alpha=\beta=1$, the
probability of seeing another sunrise after seeing $K$ out of $K$
sunrises is thus:
\[ \begin{array}{rcccl}
\mean\Big(\binomial[K]^{\dag}_{\Beta(1,1)}(K)\Big)
& = &
\displaystyle\frac{1+K}{1+1+K}
& = &
\displaystyle\frac{K+1}{K+2}.
\end{array}  \]

\noindent This is what Laplace calculated.

There is a similar situation for multinomial
distributions~\eqref{MultinomialEqn}, where the appropriate
distribution is the Dirichlet distribution $\Dirichlet(\psi)$ on
$\Dst(X)$, for a finite set $X$. This distribution involves a multiset
$\psi\in\Mlt(X)$ with full support as parameter.

\begin{lemma}
\label{DirichletLem}
\begin{enumerate}
\item \label{DirichletLemMean} $\mean\big(\Dirichlet(\psi)\big) = \flrn(\psi)$.

\item \label{DirichletLemDag} The Dirichlet distribution on $\Dst(X)$
  is conjugate prior to the multinomial channel $\multinomial[K]
  \colon \Dst(X) \chanto \Mlt[K](X)$ with parameter update induced
  by observation $\varphi\in\Mlt[K](X)$ is:
\[ \begin{array}{rcl}
\multinomial[K]^{\dag}_{\Dirichlet(\psi)}(\varphi)
& = &
\Dirichlet\big(\psi + \varphi\big).
\end{array} \]

\item \label{DirichletLemLaplace} Laplace's rule of succession thus
  gives:
\[ \begin{array}{rcccl}
\mean\Big(\multinomial[K]^{\dag}_{\Dirichlet(\psi)}(\varphi)\Big)
& = &
\mean\Big(\Dirichlet\big(\psi + \varphi\big)\Big)
& = &
\flrn\big(\psi+\varphi\big).
\end{array} \eqno{\QEDbox} \]
\end{enumerate}
\end{lemma}

We now come to the bivariate binomial situation.

\begin{proposition}
\label{BivBinDirichletProp}
For $K\in\NNO$ consider the bivariate binomial as a channel
$\bivbinomial[K] \colon \Dst\big(\finset{2}\times\finset{2}\big)
\rightarrow \{0,\ldots,K\}\times\{0,\ldots,K\}$, with a Dirichlet
distribution $\Dir(\psi)$ on
$\Dst\big(\finset{2}\times\finset{2}\big)$, for a multiset parameter
$\psi\in\Mlt\big(\finset{2}\times\finset{2}\big)$ with full
support. Laplace's rule of succession gives for observations $n_{1},
n_{2} \in \{0,\ldots,K\}$ as outcome:
\[ \begin{array}{rcl}
\mean\Big(\bivbinomial[K]^{\dag}_{\Dirichlet(\psi)}(n_{1},n_{2})\Big)
& = &
\displaystyle\flrn\left(
   \sum_{\varphi\in\Mlt[K](\finset{2}\times\finset{2}), \, \marginalheads(\varphi) = (n_{1}, n_{2})} \, \psi + \varphi\right).
\end{array} \]
\end{proposition}

\begin{myproof} (Sketch)
We write $Z$ for a suitable normalisation constant in:
\[ \begin{array}[b]{rcl}
\lefteqn{\mean\Big(\bivbinomial[K]^{\dag}_{\Dirichlet(\psi)}(n_{1},n_{2})\Big)}
\\[+0.2em]
& = &
\displaystyle \frac{1}{Z} \cdot \int_{\gamma\in\Dst(\finset{2}\times\finset{2})} 
   \bivbinomial[K](\gamma)(n_{1}, n_{2}) \cdot 
      \Dirichlet(\psi)(\gamma) \cdot \gamma \intd \gamma
\\[+1em]
& = &
\displaystyle \frac{1}{Z} \cdot 
   \sum_{\varphi\in\Mlt[K](\finset{2}\times\finset{2}), 
   \, \marginalheads(\varphi) = (n_{1}, n_{2})} \, 
   \int_{\gamma\in\Dst(\finset{2}\times\finset{2})} 
   \multinomial[K](\gamma)(\varphi) \cdot 
      \Dirichlet(\psi)(\gamma) \cdot \gamma \intd \gamma
\\[+1em]
& = &
\displaystyle\flrn\left(
   \sum_{\varphi\in\Mlt[K](\finset{2}\times\finset{2}), 
   \, \marginalheads(\varphi) = (n_{1}, n_{2})} \, 
   \mean\Big(\multinomial[K]^{\dag}_{\Dirichlet(\psi)}(\varphi)\Big)\right)
\\[+1em]
& = &   
\displaystyle\flrn\left(
   \sum_{\varphi\in\Mlt[K](\finset{2}\times\finset{2}), \, \marginalheads(\varphi) = (n_{1}, n_{2})} \, \flrn\big(\psi + \varphi\big)\right),
   \quad\mbox{by Lemma~\ref{DirichletLem}~\eqref{DirichletLemLaplace}}
\\[+1em]
& = &
\displaystyle\flrn\left(
   \sum_{\varphi\in\Mlt[K](\finset{2}\times\finset{2}), \, \marginalheads(\varphi) = (n_{1}, n_{2})} \, \psi + \varphi\right).
\end{array} \eqno{\QEDbox} \]
\end{myproof}

\ignore{

# Dirichlet experiments

def BivariateBinomial_args(K, sp):
    return AChannel(lambda *xs: BivariateBinomial(K)(DState(list(xs), sp)), 1)

def Multinomial_args(K, sp):
    return AChannel(lambda *xs: Multinomial(K)(DState(list(xs), sp)), 1)

TwoTwo = range_sp(2) @ range_sp(2)
X = Space("00", "01", "10", "11")

m = DState([1,3,4,2], X)

print( m )
mt = DState(m.matrix, TwoTwo)

K = 3
K2 = range_sp(K+1) @ range_sp(K+1)

print("\nDraw size ", K, " and multiset ", m)
print("")
#print( BivariateBinomial_args(K,TwoTwo) >> Dirichlet(m) )

# Draw size  2  and multiset  2|00> + 3|01> + 1|10> + 2|11>
# 
# 0.0833|0, 0> + 0.167|0, 1> + 0.167|0, 2> + 0.0556|1, 0> +
# 0.194|1, 1> + 0.167|1, 2> + 0.0278|2, 0> + 0.0556|2, 1> + 0.0833|2, 2>
#
# This is marginal heads after Polya
#
# 1/12|[0, 0]> + 1/6|[0, 1]> + 1/6|[0, 2]> + 1/18|[1, 0]> +
# 7/36|[1, 1]> + 1/36|[2, 0]> + 1/6|[1, 2]> + 1/18|[2, 1]> + 1/12|[2, 2]>

# m = DState([2,3,1,2], TwoTwo)
# print( Functor(lambda m: [marginal_heads(m)[0], marginal_heads(m)[1]])
#        (Polya(K)(m)) )
# print("")
# print( Functor(lambda m: [marginal_heads(m)[0], marginal_heads(m)[1]])
#        (Polya(K, frac=True)(m)) )

#
#
#
n1 = 3
n2 = 1
p = AFunction(lambda *xs: BivariateBinomial(K)(DState(list(xs), TwoTwo))(n1,n2), 4)
#print( p )
# print( Dirichlet(m) / p >= AFunction(lambda xs: xs[0], 4) )
# print( Dirichlet(m) / p >= AFunction(lambda xs: xs[1], 4) )
# print( Dirichlet(m) / p >= AFunction(lambda xs: xs[2], 4) )
# print( Dirichlet(m) / p >= AFunction(lambda xs: xs[3], 4) )

#
# In accordance with the listing below
#
# 2|00> + 3|01> + 1|10> + 2|11>,  K=2, n1=2, n2=1
# 0.20000000095045356
# 0.300000001728039
# 0.20000000107651994
# 0.30000000152219186

# Draw size  3  and multiset  1|00> + 3|01> + 4|10> + 2|11>
# at 00: 0.07692307752494579, K=3, n1=3, n2=1, corresponds with below

# 
for n1 in range(K+1):
    for n2 in range(K+1):
        print(n1, n2, "  ", list_addition(
            [ mt + m1 for m1 in Multisets(K)(TwoTwo)
              if marginal_heads(m1) == (n1, n2) ]).flrn() )

}

In Lemma~\ref{BetaLem} and~\ref{DirichletLem} there are Laplace
succession rules arising from conjugate prior situations, but this not
a necessary ingredient, as Equation~\eqref{PoissonLaplaceEqn} below
shows.  It involves the combination of the binomial and Poisson
distributions, see \textit{e.g.}~\cite[\S6.11]{Jaynes03}. This applies
for instance when one has an imperfect particle counter/detector in
physics, where particles are emitted by some source according to a
Poisson distribution (with rate $\lambda$) and the imperfect detection
is captured by a subsequent binomial channel (with chance $r$ of being
detected). When $n$ particles are detected, the question is: what is
the expected number of emitted particles?  The distribution of these
total numbers is obtained by updating the prior Poisson
distribution. Interestingly, this yields a new Poisson distribution,
but with a shift. The Laplace succession rule says in this case:
\begin{equation}
\label{PoissonLaplaceEqn}
\begin{array}{rcl}
\mean\Big(\binomial[-](r)^{\dag}_{\poisson[\lambda]}(n)\Big)
& = &
n + (1-r)\cdot\lambda.
\end{array}
\end{equation}

\noindent Here one views binomial as a channel $\binomial[-](r) \colon
\NNO \chanto \{0,\ldots,K\}$ that is reversed via the dagger, with the
Poisson distribution on $\NNO$ as prior. The
formula~\eqref{PoissonLaplaceEqn} expresses that after detecting $n$
particles the expected number of emitted particles is $n +
(1-r)\cdot\lambda$. The factor $(1-r)\cdot\lambda$ captures the
expected number of emitted particles that were not detected. This
makes sense.

Without proof we mention the following analogue for 
bivariate binomials.

\begin{proposition}
\label{BivBinPoissonProp}
Fix $\lambda\in\nnR$, $K\in\NNO$,
$\gamma\in\Dst\big(\finset{2}\times\finset{2}\big)$. For observations
$n_{1}, n_{2} \in \{0,\ldots,K\}$, say with $n_{1} \leq n_{2}$, one has:
\[ \begin{array}{rcl}
\mean\Big(\bivbinomial[-](\gamma)^{\dag}_{\poisson[\lambda]}(n_{1}, n_{2})\Big)
& = &
\displaystyle n_{2} \,+\, \gamma(0,0)\cdot\lambda \,+\,
   \frac{\sum_{0\leq i\leq n_{1}} \mathsl{ppp}(i)\cdot i}
   {\sum_{0\leq i\leq n_{1}} \mathsl{ppp}(i)},
\end{array} \]

\noindent where $\mathsl{ppp}(i) \coloneqq
\poisson[\gamma(0,1)\cdot\lambda](n_{2}-n_{1}+i) \cdot
\poisson[\gamma(1,0)\cdot\lambda](i) \cdot
\poisson[\gamma(1,1)\cdot\lambda](n_{1}-i)$. \QED
\end{proposition}

\auxproof{

\begin{myproof}
Let's write $\mathsl{ppp} \coloneqq \sum_{0\leq i\leq n_{1}}
\mathsl{ppp}(i)$.  Our first goal is to prove that
$\pushing{\bivbinomial[-](\gamma)}{\poisson[\lambda]}(n_{1}, n_{2}) =
\mathsl{ppp}$, where we still assume $n_{1} \leq n_{2}$. We
use Equation~\eqref{MarginalHeadsInverseEqn} in:
\[ \begin{array}{rcl}
\lefteqn{\pushing{\bivbinomial[-](\gamma)}{\poisson[\lambda]}(n_{1}, n_{2})}
\\[+0.2em]
& = &
\displaystyle\sum_{K\geq n_{2}} \, \bivbinomial[K](\gamma)(n_{1}, n_{2}) 
   \cdot\poisson[\lambda](K)
\\[+1em]
& = &
\displaystyle\sum_{K\geq n_{2}} \;
   \sum_{\varphi\in\Mlt[K](\finset{2}\times\finset{2}), \,
   \marginalheads(\varphi) = (n_{1},n_{2})} \,
   \multinomial[K](\gamma)(\varphi) \cdot\poisson[\lambda](K)
\\[+1em]
& = &
\displaystyle\sum_{K\geq n_{2}} \; \sum_{0 \leq i \leq \min(n_{1}, K-n_{2})} \,
   \frac{K!}{(K - n_{2} - i)!\cdot (n_{2} - n_{1} + i)! 
   \cdot i! \cdot (n_{1} - i)!}
\\
& & \qquad \displaystyle \cdot 
   \gamma(0,0)^{K - n_{2} - i} \cdot \gamma(0,1)^{n_{2} - n_{1} + i} \cdot
   \gamma(1,0)^{i} \cdot \gamma(1,1)^{n_{1} - i} \cdot
   e^{-\lambda} \cdot \frac{\lambda^{K}}{K!}
\\[+1em]
& = &
\displaystyle\sum_{K\geq n_{2}} \; \sum_{0 \leq i \leq \min(n_{1}, K-n_{2})} \,
   e^{-\gamma(0,0)\cdot\lambda} \cdot 
   \frac{(\gamma(0,0)\cdot\lambda)^{K - n_{2} - i}}{(K - n_{2} - i)!} \cdot
   e^{-\gamma(0,1)\cdot\lambda} \cdot 
   \frac{(\gamma(0,1)\cdot\lambda)^{n_{2} - n_{1} + i}}{(n_{2} - n_{1} + i)!}
\\
& & \qquad \displaystyle \cdot 
   e^{-\gamma(1,0)\cdot\lambda} \cdot 
   \frac{(\gamma(1,0)\cdot\lambda)^{i}}{i!} \cdot
   e^{-\gamma(1,1)\cdot\lambda} \cdot 
   \frac{(\gamma(1,1)\cdot\lambda)^{n_{1} - i}}{(n_{1} - i)!}
\\[+1em]
& = &
\displaystyle\sum_{K\geq n_{2}} \; \sum_{0 \leq i \leq \min(n_{1}, K-n_{2})} \,
   \poisson[\gamma(0,0)\cdot\lambda](K - n_{2} - i) \cdot \mathsl{ppp}(i)
\\[+1em]
& = &
\displaystyle\sum_{n_{2} \leq K \leq n_{1} + n_{2}} \; 
   \sum_{0 \leq i \leq K-n_{2}} \,
   \poisson[\gamma(0,0)\cdot\lambda](K - n_{2} - i) \cdot \mathsl{ppp}(i)
\\
& & \qquad \displaystyle + \,
\displaystyle\sum_{K > n_{1} + n_{2}} \; 
   \sum_{0 \leq i \leq n_{1}} \,
   \poisson[\gamma(0,0)\cdot\lambda](K - n_{2} - i) \cdot \mathsl{ppp}(i)
\\[+1em]
& = &
\poisson[\gamma(0,0)\cdot\lambda](0) \cdot \mathsl{ppp}(0)
\\
& & \qquad + \;
\poisson[\gamma(0,0)\cdot\lambda](1) \cdot \mathsl{ppp}(0) + 
\poisson[\gamma(0,0)\cdot\lambda](0) \cdot \mathsl{ppp}(1)
\\
& & \qquad +\; \cdots 
\\
& & \qquad + \;
\poisson[\gamma(0,0)\cdot\lambda](n_{1}) \cdot \mathsl{ppp}(0) + \cdots +
\poisson[\gamma(0,0)\cdot\lambda](0) \cdot \mathsl{ppp}(n_{1})
\\
& & \qquad +\; \displaystyle\sum_{K > n_{1}} \; 
   \sum_{0 \leq i \leq n_{1}} \,
   \poisson[\gamma(0,0)\cdot\lambda](K - i) \cdot \mathsl{ppp}(i)
\\[+1em]
& = &
\displaystyle\sum_{K\geq 0} \; \sum_{0\leq i \leq n_{1}} \,
   \poisson[\gamma(0,0)\cdot\lambda](K) \cdot \mathsl{ppp}(i)
\\
& = &
\mathsl{ppp}.
\end{array} \]

We can now come to the mean of the update:
\[ \begin{array}{rcl}
\lefteqn{\mean\Big(\bivbinomial[-](\gamma)^{\dag}_{\poisson[\lambda]}(n_{1}, n_{2})\Big)}
\\[+0.2em]
& = &
\displaystyle\sum_{K\geq n_{2}} \, 
   \frac{\bivbinomial[K](\gamma)(n_{1}, n_{2}) 
   \cdot\poisson[\lambda](K) \cdot K}
  {\pushing{\bivbinomial[-](\gamma)}{\poisson[\lambda]}(n_{1}, n_{2})}
\\[+1em]
& = &
\displaystyle\sum_{K\geq n_{2}} \;
   \sum_{\varphi\in\Mlt[K](\finset{2}\times\finset{2}), \,
   \marginalheads(\varphi) = (n_{1},n_{2})} \,
   \frac{\multinomial[K](\gamma)(\varphi) \cdot\poisson[\lambda](K)\cdot K}
   {\mathsl{ppp}}
\\[+1em]
& = &
\displaystyle\sum_{K\geq n_{2}} \; \sum_{0 \leq i \leq \min(n_{1}, K-n_{2})} \,
   \frac{K!}{(K - n_{2} - i)!\cdot (n_{2} - n_{1} + i)! 
   \cdot i! \cdot (n_{1} - i)!}
\\
& & \qquad \displaystyle \cdot 
   \frac{\gamma(0,0)^{K - n_{2} - i} \cdot \gamma(0,1)^{n_{2} - n_{1} + i} \cdot
   \gamma(1,0)^{i} \cdot \gamma(1,1)^{n_{1} - i}}{\mathsl{ppp}} \cdot
   e^{-\lambda} \cdot \frac{\lambda^{K}}{K!} \cdot K
\\[+1em]
& = &
\displaystyle\sum_{K\geq n_{2}} \; \sum_{0 \leq i \leq \min(n_{1}, K-n_{2})} \,
   \frac{\poisson[\gamma(0,0)\cdot\lambda](K - n_{2} - i) \cdot 
   \mathsl{ppp}(i)\cdot K}{\mathsl{ppp}}
\\[+1em]
& = &
\displaystyle \sum_{K \geq \NNO} \; \sum_{0\leq i\leq n_{1}} \,
   \frac{\poisson[\gamma(0,0)\cdot\lambda](K) \cdot 
   \mathsl{ppp}(i)\cdot (K + n_{2} + i)}{\mathsl{ppp}}
\\[+1em]
& = &
\displaystyle \sum_{K \geq \NNO} \; \sum_{0\leq i\leq n_{1}} \,
   \frac{\poisson[\gamma(0,0)\cdot\lambda](K) \cdot 
   \mathsl{ppp}(i)\cdot K}{\mathsl{ppp}} \;+\;
   \sum_{K \geq \NNO} \; \sum_{0\leq i\leq n_{1}} \,
   \frac{\poisson[\gamma(0,0)\cdot\lambda](K) \cdot 
   \mathsl{ppp}(i)\cdot n_{2}}{\mathsl{ppp}}
\\
& & \qquad \displaystyle + \,
\sum_{K \geq \NNO} \; \sum_{0\leq i\leq n_{1}} \,
   \frac{\poisson[\gamma(0,0)\cdot\lambda](K) \cdot 
   \mathsl{ppp}(i)\cdot i}{\mathsl{ppp}}
\\[+1em]
& = &
\mean\Big(\poisson[\gamma(0,0)\cdot\lambda]\Big) \;+\; n_{2} \;+\; 
   \displaystyle \sum_{0\leq i\leq n_{1}} \, 
   \frac{\mathsl{ppp}(i)\cdot i}{\mathsl{ppp}}
\\[+1em]
& = &
\displaystyle n_{2} \,+\, \gamma(0,0)\cdot\lambda \,+\,
   \frac{\sum_{0\leq i\leq n_{1}} \mathsl{ppp}(i)\cdot i}
   {\sum_{0\leq i\leq n_{1}} \mathsl{ppp}(i)}.
\end{array} \]
\end{myproof}

}

\ignore{

print("\nBivariate binomial and Poisson")

TwoTwo = range_sp(2) @ range_sp(2)
n2 = random.randint(0,10)
n1 = random.randint(0,n2)
print("\nHeads: ", n1, n2)

w = random_distribution(TwoTwo)
lam = random.uniform(0,5)
M = 15
pl = Poisson(lam, MaxRange=M)

def ppp(i):
    return Poisson(w(1,0) * lam)(i) * \
        Poisson(w(0,1) * lam)(n2-n1+i) * \
        Poisson(w(1,1) * lam)(n1-i)

pppf = sum([ ppp(i) for i in range(0,n1+1) ])

print( pppf )

print( sum([ sum([ Poisson(w(0,0) * lam)(K-n2-i) *
                   Poisson(w(1,0) * lam)(i) *
                   Poisson(w(0,1) * lam)(n2-n1+i) *
                   Poisson(w(1,1) * lam)(n1-i)
                   for i in range(min(n1, K-n2)+1) ])
             for K in range(n2, M+1) ]) )

print( sum([ sum([ Poisson(w(0,0) * lam)(K-n2-i) *
                   Poisson(w(1,0) * lam)(i) *
                   Poisson(w(0,1) * lam)(n2-n1+i) *
                   Poisson(w(1,1) * lam)(n1-i)
                   for i in range(K-n2+1) ])
             for K in range(n2,n1+n2) ])
       + sum([ sum([ Poisson(w(0,0) * lam)(K-n2-i) *
                   Poisson(w(1,0) * lam)(i) *
                   Poisson(w(0,1) * lam)(n2-n1+i) *
                   Poisson(w(1,1) * lam)(n1-i)
                   for i in range(n1+1)  ])
               for K in range(n1+n2,M+1) ]) )

# Somewhat deviant outcome
print( sum([ sum([ Multinomial(K)(w)(m)
             for m in Multisets(K)(TwoTwo)
             if marginal_heads(m) == (n1,n2) ]) * Poisson(lam)(K)
             for K in range(n2, M+1) ]) )

# commented out, since slow

print( (AChannel(lambda K: Functor(lambda x: x, cod=range_sp(M+1)**2)
                 (BivariateBinomial(K)(w)),
                 pl.sp,
                 range_sp(M+1)**2) >> pl)(n1, n2) )

p = DPred(lambda K: BivariateBinomial(K)(w)(n1,n2) if K >= n2 else 0, pl.sp)
print( "\nMean of update")
e = (pl / p).expectation()
print( e )
# same as:
print( sum([ pl(k) * BivariateBinomial(k)(w)(n1,n2) * k
             for k in range(n2, M+1) ]) / pppf )
print( n2 + w(0,0) * lam + sum([ ppp(i) * i
                                 for i in range(n1+1) ]) / pppf )

}

\section{Expectation Maximisation for bivariate binomials}\label{EMSec}

Expectation Maximisation is an unsupervised machine learning technique
for classifying data in a finite number of classes/features. The
classification takes place by finding a `mixture' distribution that
assigns propertions (in the form of probabilities) to the different
classification features.  This section gives an impression how
bivariate distributions --- in the formulation of this paper --- may
be useful in concrete applications in machine learning. Below we use
the channel-based analysis of Expectation Maximisation developed
in~\cite{Jacobs23b,Jacobs21g} and demonstrate how it works in a
concrete example. It avoids explicit updates of parameters, as happens
in common formulations in Expectation Maximisation.  However,
implicitly, a two-coin distribution parameter is recomputed in each
iteration via the steps described in Fact~\ref{RecoverFact}.

\begin{figure}
\begin{center}
  \includegraphics[scale=0.4]{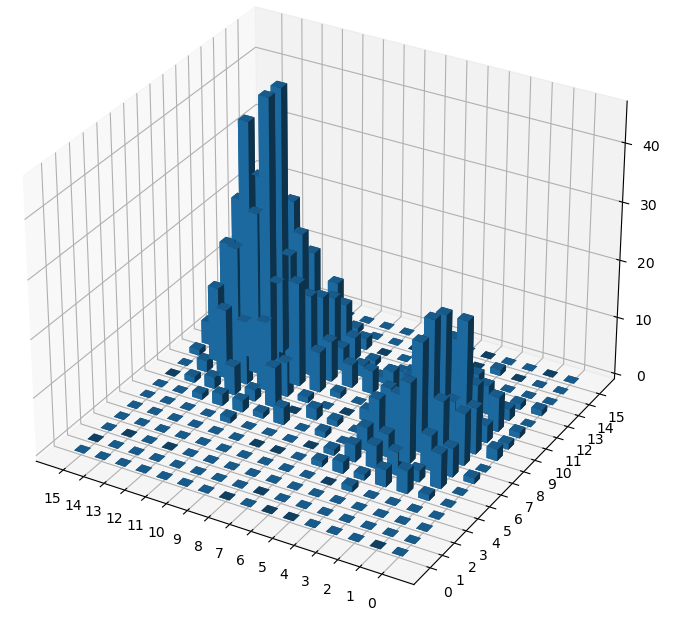}
\end{center}
\caption{Plots of 1000 samples from the distribution $\sigma\in
  \Dst\big(\{0,\ldots,15\}\times\{0,\ldots,15\}\big)$
  from~\eqref{MixtureEqn}.}
\label{DataFig}
\end{figure}

The approach is as follows. We start from a known mixture of two
bivariate binomial distributions, with $K=15$, namely:
\begin{equation}
\label{MixtureEqn}
\hspace*{-0.5em}\begin{array}{rcl}
\sigma
& \coloneqq &
\frac{1}{3} \cdot \bivbinomial[K](\gamma_{0}) \,+\,
   \frac{2}{3} \cdot \bivbinomial[K](\gamma_{1})
\;\mbox{ where }\;
\left\{\begin{array}{rcl}
\gamma_{0} 
& = &
\frac{3}{8}\bigket{0,0} +
   \frac{5}{12}\bigket{0,1} +
   \frac{1}{12}\bigket{1,0} +
   \frac{1}{8}\bigket{1,1}
\\[+0.2em]
\gamma_{1} 
& = &
\frac{1}{10}\bigket{0,0} +
   \frac{1}{10}\bigket{0,1} +
   \frac{1}{5}\bigket{1,0} +
   \frac{3}{5}\bigket{1,1}
\end{array}\right.
\end{array}
\end{equation}

\noindent We form a multiset $\psi \in
\Mlt[1000]\big(\finset{2}\times\finset{2}\big)$ via 1000 random
samples from this distribution. This multiset is plotted in
Figure~\ref{DataFig}. One recognises two humps that indicate that
there may be an underlying mixture of bivariate binomial
distributions. Our aim is now to reconstruct not only the $\frac{1}{3}
- \frac{2}{3}$ mixture but also the two distributions $\gamma_{1},
\gamma_{2}$ from these data alone --- where we only assume that two
numbers are known, namely the number $K=15$ of tosses and the
number~$2$ of binomials in the mixture.

\begin{figure}
\[ \vcenter{\hbox{\begin{minipage}[t]{0.88\textwidth}
\begin{lstlisting}
def BivariateBinomialEM (dist, chan, data_dist):
   dagger = chan$^{\dag}_{\text{dist}}$
   # E-step, as Jeffrey update
   new_dist = dagger$_{*}$(data_dist)
   # M-step, via double dagger
   double_dagger = dagger$^{\dag}_{\text{data\_dist}}$
   def new_chan(x):  return bvbn[K]( recover(double_dagger(x)) )
   return (new_dist, new_chan)
\end{lstlisting}
\end{minipage}}} \]

\caption{Expectation Maximisation for finding a mixture of 
bivariate binomial distributions.}
\label{ProgramFig}
\end{figure}

Expectation Maximisation is the technique for doing this. It seeks a
mixture distribution $\omega$, in this case on $\finset{2} = \{0,1\}$
since we seek a mixture of~$2$ binomials, and a channel $c\colon
\finset{2} \chanto \{0,\ldots,15\}^{2}$, where both $c(0)$ and $c(1)$
are bivariate distributions. The aim is to minimise the
KL-divergence $\DKL\big(\flrn(\psi), \pushing{c}{\omega}\big)$, so
that the divergence between the data distribution and the prediction
is minimal, see~\cite{Jacobs23b} for details.

This reduction happens via several iterations. A single iteration, in
Python-like pseudo code, is described in Figure~\ref{ProgramFig}.  It
produces a new mixture distribution via Jeffrey's update
rule~\cite{Jacobs19c,Jacobs21c}, using the dagger of the original
channel, with the mixture distribution as prior. The new channel is
obtained by first forming the double dagger, with the data
distribution $\flrn(\psi)$ as prior. This double dagger may not be of
the right bivariate binomial form, but it is forced into this
form via a \inline{recover} function. It uses the steps described in
Fact~\ref{RecoverFact}, to find a distribution on
$\finset{2}\times\finset{2}$ by computing the mean and covariance.

We run the algorithm in Figure~\ref{ProgramFig} five times on an
arbitrary mixture distribution $\omega\in\Dst(\finset{2})$ and with a
channel $c\colon \finset{2} \chanto \{0,\ldots,15\}^{2}$ consisting of
arbitrary bivariate distributions $c(0)$ and $c(1)$. One such
run gives as successive divergences:
\[ 
1.991
\qquad
0.272
\qquad
0.093
\qquad
0.087
\qquad
0.087.
\] 

\noindent Starting from a divergence resulting from random choices, we
see that these divergence quickly stabilise (in three decimals). The
mixture distribution and associated (channel of bivariate
binomials) given by two-coin distributions that emerge after these
five runs are:
\[ \begin{array}{c}
0.326\bigket{0} + 0.674\bigket{1}
\end{array}
\hspace*{4em}
\begin{array}{c}
0.378\bigket{0, 0} + 0.416\bigket{0, 1} + 
   0.0795\bigket{1, 0} + 0.127\bigket{1, 1}
\\
0.094\bigket{0, 0} + 0.104\bigket{0, 1} + 
   0.205\bigket{1, 0} + 0.597\bigket{1, 1}
\end{array} \]

\noindent We see that they are pretty close to values in the
mixture~\eqref{MixtureEqn}, from which the 1000 samples in
Figure~\ref{DataFig} were taken. This example illustrates how
Expectation Maximisation can be performed for bivariate binomial
distributions.

\ignore{

# Expectation Maximisation experiment

TwoTwo = range_sp(2) ** 2
K = 10
w = random_distribution(TwoTwo)
print( w )
bt = BivariateBinomial(K)(w)
print( recover_twocoin_from_bivariate_binomial(bt, K) )

#BT.dplot2()

#psi = BT.sample(1000)
#print( psi.matrix )

# [[ 0  0  0  0  0  0  1  0  2  1  1  0  1  0  0  0]
#  [ 0  0  0  0  1  6  5  9  8  3  6  2  1  0  0  0]
#  [ 0  0  0  0  4  2  6 10  8 20  7  5  1  0  0  0]
#  [ 0  0  0  0  3  4 14 19 21 20 14  7  4  0  1  0]
#  [ 0  0  0  1  1  4  4  9 10 11  8  3  2  1  0  0]
#  [ 0  0  0  0  2  3  4  7  7  3  6  1  1  0  1  0]
#  [ 0  0  0  0  1  1  0  0  2  2  5  3  0  0  0  0]
#  [ 0  0  0  0  0  0  1  1  0  0  4  1  0  0  0  0]
#  [ 0  0  0  0  0  0  0  2  0  1  4  2  1  0  0  0]
#  [ 0  0  0  0  0  0  3  0  1  7  7  4  4  2  0  0]
#  [ 0  0  0  0  0  0  1  7  6 18 14 12 10  7  1  0]
#  [ 0  0  0  0  0  1  2  2 12 17 20 22 17  5  8  0]
#  [ 0  0  0  0  0  0  2  5 11 28 46 46 25 14  3  1]
#  [ 0  0  0  0  0  0  1  2 12 21 41 30 33 10  3  0]
#  [ 0  0  0  0  0  0  0  1  2 13 19 25 16 13  5  2]
#  [ 0  0  0  0  0  0  0  0  0  1  4  8  6  3  1  0]]

data = DState([
    [ 0,  0,  0,  0,  0,  0,  1,  0,  2,  1,  1,  0,  1,  0,  0,  0],
    [ 0,  0,  0,  0,  1,  6,  5,  9,  8,  3,  6,  2,  1,  0,  0,  0],
    [ 0,  0,  0,  0,  4,  2,  6, 10,  8, 20,  7,  5,  1,  0,  0,  0],
    [ 0,  0,  0,  0,  3,  4, 14, 19, 21, 20, 14,  7,  4,  0,  1,  0],
    [ 0,  0,  0,  1,  1,  4,  4,  9, 10, 11,  8,  3,  2,  1,  0,  0],
    [ 0,  0,  0,  0,  2,  3,  4,  7,  7,  3,  6,  1,  1,  0,  1,  0],
    [ 0,  0,  0,  0,  1,  1,  0,  0,  2,  2,  5,  3,  0,  0,  0,  0],
    [ 0,  0,  0,  0,  0,  0,  1,  1,  0,  0,  4,  1,  0,  0,  0,  0],
    [ 0,  0,  0,  0,  0,  0,  0,  2,  0,  1,  4,  2,  1,  0,  0,  0],
    [ 0,  0,  0,  0,  0,  0,  3,  0,  1,  7,  7,  4,  4,  2,  0,  0],
    [ 0,  0,  0,  0,  0,  0,  1,  7,  6, 18, 14, 12, 10,  7,  1,  0],
    [ 0,  0,  0,  0,  0,  1,  2,  2, 12, 17, 20, 22, 17,  5,  8,  0],
    [ 0,  0,  0,  0,  0,  0,  2,  5, 11, 28, 46, 46, 25, 14,  3,  1],
    [ 0,  0,  0,  0,  0,  0,  1,  2, 12, 21, 41, 30, 33, 10,  3,  0],
    [ 0,  0,  0,  0,  0,  0,  0,  1,  2, 13, 19, 25, 16, 13,  5,  2],
    [ 0,  0,  0,  0,  0,  0,  0,  0,  0,  1,  4,  8,  6,  3,  1,  0]],
              range_sp(16) @ range_sp(16))

data_dist = data.flrn()

#print( data )
#data.dplot2()

K = 15
w = DState([3/8, 5/12, 1/12, 1/8], TwoTwo)
v = DState([1/10, 1/10, 1/5, 3/5], TwoTwo)
#BT = 1/3 * BivariateBinomial(K)(w) + 2/3 * BivariateBinomial(K)(v)

print("\nUsed in 1/3 - 2/3 mixture")
print(" * ", w )
print(" * ", v )
print("")

mixt = random_distribution(range_sp(2))
chan = DChannel([ BivariateBinomial(K)(random_distribution(TwoTwo)),
                  BivariateBinomial(K)(random_distribution(TwoTwo)) ],
                range_sp(2), range_sp(K+1) ** 2)

for i in range(5):
    mixt, chan = BivariateBinomialEM(mixt, chan, K, data_dist)

print("\nResults:")
print( mixt )
print( recover_twocoin_from_bivariate_binomial(chan(0), K) )
print( recover_twocoin_from_bivariate_binomial(chan(1), K) )
print("Final divergence: ", round(kldive(data_dist, chan >> mixt), 3) )

# Divergence:  1.991
# Divergence:  0.272
# Divergence:  0.093
# Divergence:  0.087
# Divergence:  0.087

# Results:
# 0.326|0> + 0.674|1>
# 0.378|0, 0> + 0.416|0, 1> + 0.0795|1, 0> + 0.127|1, 1>
# 0.094|0, 0> + 0.104|0, 1> + 0.205|1, 0> + 0.597|1, 1>
# Final divergence:  0.087

}

\section{Concluding remarks}

This paper reformulates bivariate binomial distributions in a concise
manner, by exploiting the functoriality of the mapping $X \mapsto
\Dst(X)$, sending a set $X$ to the set of probability distributions on
$X$.  It demonstrates the clarifying power of the categorical approach
to probability theory by (re)describing some of the basic properties
of bivariate binomials, for instance via (dagger) channels. It also
shows that this categorical formulation makes it possible to prove
basic properties via equational reasoning. The paper concentrates on
two dimensions, for simplicity, but the extension to multiple
dimensions is straightforward.

The multivariate binomial distribution be understood as a discrete
version of the widely used multidimensional Gaussian distribution.  By
the famous theorem of De Moivre-Laplace every binomial distributions
approximate a Gaussian distribution. This is a special case of the
Central Limit Theorem. The approximation also works for multivariate
binomials, see~\cite{Veeh86}.








\appendix

\section*{Appendix}

We add proofs that are skipped in the main text.

\begin{myproof} (of Lemma~\ref{BivBinMarginalLem})
\begin{enumerate}
\item We do only the first equation, using twice the Binomial
  Theorem. For $0\leq n_{1}\leq K$ we have:
\[ \begin{array}{rcl}
\Dst\big(\Pi_{1}\big)\Big(\bivbinomial[K](\gamma)\Big)(n_{1})
& = &
\displaystyle\sum_{0\leq n_{2}\leq K} \,
   \Dst\big(\marginalheads\big)\Big(\multinomial[K](\gamma)\Big)(n_{1}, n_{2})
\\[+1em]
& = &
\displaystyle\sum_{0\leq n_{2}\leq K} \;
   \sum_{\varphi\in\Mlt[K](\finset{2}\times \finset{2}), \, 
   \marginalheads(\varphi) = (n_{1}, n_{2})} \,
   \multinomial[K](\gamma)(\varphi)
\\[+1.4em]
& = &
\displaystyle \sum_{\varphi\in\Mlt[K](\finset{2}\times \finset{2}), \, 
   \varphi(1,0) + \varphi(1,1) = n_{1}} \,
   \coefm{\varphi} \cdot \prod_{i,j\in\finset{2}} \, \gamma(i,j)^{\varphi(i,j)}
\\[+1.2em]
& = &
\displaystyle \sum_{0 \leq i \leq n_{1}} \, \sum_{0 \leq j \leq K-n_{1}} \,
    \frac{K!}{(K - n_{1} - j)! \cdot j! \cdot (n_{1}-i)! \cdot i!} 
\\
& & \qquad \cdot\,
\gamma(0,0)^{K - n_{1} - j} \cdot \gamma(0,1)^{j} \cdot
   \gamma(1,0)^{n_{1}-i} \cdot \gamma(1,1)^{i}
\\
& = &
\displaystyle \sum_{0 \leq i \leq n_{1}} \, \sum_{0 \leq j \leq K-n_{1}} \,
    \frac{K!}{n_{1}! \cdot (K-n_{1})!} \cdot
    \frac{n_{1}!}{i!\cdot (n_{1}-i)!} \cdot 
    \frac{(K-n_{1})!}{j! \cdot (K-n_{1}-j)!}
\\
& & \qquad \cdot\,
\gamma(0,0)^{K - n_{1} - j} \cdot \gamma(0,1)^{j} \cdot
   \gamma(1,0)^{n_{1}-i} \cdot \gamma(1,1)^{i}
\\[+0.2em]
& = &
\displaystyle \binom{K}{n_1} \cdot \big(\gamma(1,0) + \gamma(1,1)\big)^{n_1} 
   \cdot \big(\gamma(0,0) + \gamma(0,1)\big)^{K-n_1} 
\\[+0.8em]
& = &
\displaystyle \binom{K}{n_1} \cdot \gamma_{1}(1)^{n_1} \cdot
   \big(1 - \gamma_{1}(1)\big)^{K-n_1}
\hspace*{\arraycolsep}=\hspace*{\arraycolsep}
\binomial[K]\Big(\gamma_{1}(1)\Big)(n_{1}).
\end{array} \]

\item We now get for $0\leq n_{1}, n_{2}\leq K$, via
  Lemma~\ref{MarginalHeadsCoefficientLem},
\[ \hspace*{-1em}\begin{array}[b]{rcl}
\bivbinomial[K]\Big(\gamma_{1}\otimes\gamma_{2}\Big)(n_{1},n_{2})
& = &
\displaystyle\sum_{\varphi\in\Mlt[K](\finset{2}\times \finset{2}), \, 
   \marginalheads(\varphi) = (n_{1}, n_{2})} \, 
   \multinomial[K]\Big(\gamma_{1}\otimes\gamma_{2}\Big)(\varphi)
\\[+1.2em]
& = &
\displaystyle\sum_{\varphi\in\marginalheads^{-1}(n_{1},n_{2})} \,
\coefm{\varphi} \cdot
\gamma_{1}(1)^{\varphi(1,0)+\varphi(1,1)} \cdot
\gamma_{1}(0)^{\varphi(0,0) + \varphi(0,1)} \cdot
\\
& & \qquad \cdot\,
\gamma_{2}(1)^{\varphi(0,1)+\varphi(1,1)} \cdot
\gamma_{2}(0)^{\varphi(0,0) + \varphi(1,0)}
\\[+0.2em]
& = &
\displaystyle\sum_{\varphi\in\marginalheads^{-1}(n_{1},n_{2})} \,
\coefm{\varphi} \cdot
\gamma_{1}(1)^{n_{1}} \cdot \gamma_{1}(0)^{K-n_{1}} \cdot
\gamma_{2}(1)^{n_{2}} \cdot \gamma_{2}(0)^{K-n_{2}}
\\[+1.2em]
& = &
\displaystyle\binom{K}{n_{1}}\cdot\binom{K}{n_{2}} \cdot
\gamma_{1}(1)^{n_{1}} \cdot \big(1 - \gamma_{1}(1)\big)^{K-n_{1}} \cdot
\gamma_{2}(1)^{n_{2}} \cdot \big(1 - \gamma_{2}(1)\big)^{K-n_{2}}
\\[+0.2em]
& = &
\binomial[K]\big(\gamma_{1}(1)\big)(n_{1}) \cdot
   \binomial[K]\big(\gamma_{2}(1)\big)(n_{2})
\\[+0.2em]
& = &
\Big(\binomial[K]\big(\gamma_{1}(1)\big) \otimes 
   \binomial[K]\big(\gamma_{2}(1)\big)\Big)(n_{1}, n_{2}).
\end{array} \eqno{\QEDbox} \]
\end{enumerate}
\end{myproof}

\begin{myproof} (of Proposition~\ref{BivBinConvolutionProp})
For the first equation we use the closure of multinomial distributions
under convolution, from
Lemma~\ref{MultinomialLem}~\eqref{MultinomialLemConv} and note that
the commutative monoid involved is $\NNO^{2}$ with the component-wise
sum $+\cdot\NNO^{2} \times \NNO^{2} \rightarrow \NNO^{2}$.
    \[ \begin{array}{rcl}
      \lefteqn{\bivbinomial[K](\gamma) + \bivbinomial[L](\gamma)}
      \\[+0.2em]
      & = &
      \Dst(+)\Big(\bivbinomial[K](\gamma) + \bivbinomial[L](\gamma)\Big)
      \\[+0.2em]
      & = &
      \displaystyle\sum_{n_{1}, n_{2} \in \NNO} \, \sum_{m_{1},m_{2}\in\NNO} \,
      \bivbinomial[K](\gamma)(n_{1},n_{2}) \cdot
      \bivbinomial[L](\gamma)\Big)(m_{1}, m_{2})
      \Bigket{(n_{1},n_{2}) + (m_{1},m_{2})}
      \\[+0.2em]
      & = &
      \displaystyle\sum_{n_{1}, n_{2} \in \NNO} \, \sum_{m_{1},m_{2}\in\NNO} \,
      \left(\sum_{\varphi\in\marginalheads^{-1}(n_{1},n_{2})} 
      \multinomial[K](\gamma)(\varphi)\right) \cdot
      \left(\sum_{\psi\in\marginalheads^{-1}(m_{1},m_{2})} 
      \multinomial[L](\gamma)(\psi)\right)
      \\[+1.0em]
      & & \hspace*{4em}  \Bigket{(n_{1}+m_{1},n_{2}+m_{2})}
      \\[+0.2em]
      & = &
      \displaystyle\sum_{\varphi\in\Mlt[K](\finset{2}\times\finset{2})} \,
      \sum_{\psi\in\Mlt[L](\finset{2}\times\finset{2})} \,
      \multinomial[K](\gamma)(\varphi) \cdot \multinomial[L](\gamma)(\psi)
      \\[+0.8em]
      & & \hspace*{4em}  \Bigket{\big(
              \Mlt(\pi_{1})(\varphi)(1)+\Mlt(\pi_{1})(\psi)(1),
              \Mlt(\pi_{2})(\varphi)(1)+\Mlt(\pi_{2})(\psi)(1)\big)}
      \\[+0.2em]
      & = &
      \displaystyle\sum_{\varphi\in\Mlt[K](\finset{2}\times\finset{2})} \,
      \sum_{\psi\in\Mlt[L](\finset{2}\times\finset{2})} \,
      \multinomial[K](\gamma)(\varphi) \cdot \multinomial[L](\gamma)(\psi)
      \\[+1.0em]
      & & \hspace*{4em}  \Bigket{\big(
              \Mlt(\pi_{1})(\varphi+\pi_{1})(\psi)(1),
              \Mlt(\pi_{2})(\varphi+\psi)(1)\big)}
      \\[+0.4em]
      & = &
      \displaystyle \Dst\big(\marginalheads\big)\left(
      \displaystyle\sum_{\varphi\in\Mlt[K](\finset{2}\times\finset{2})} \,
      \sum_{\psi\in\Mlt[L](\finset{2}\times\finset{2})} \,
      \multinomial[K](\gamma)(\varphi) \cdot \multinomial[L](\gamma)(\psi)
      \,\Bigket{\varphi+\psi}\right)
      \\[+1.2em]
      & = &
      \Dst\big(\marginalheads\big)\Big(
      \multinomial[K](\gamma) + \multinomial[L](\gamma)\Big)
      \hspace*{\arraycolsep}=\hspace*{\arraycolsep}
      \Dst\big(\marginalheads\big)\Big(\multinomial[K+ L](\gamma)\Big)
      \hspace*{\arraycolsep}=\hspace*{\arraycolsep}
      \bivbinomial[K+ L](\gamma).
    \end{array} \]

The second equation $\bivbinomial[K](\gamma) = K\cdot \gamma$ is
obtained via this closure under convolution, using that
$\bivbinomial[0](\gamma) = 1\bigket{0,0} = 0\cdot\gamma$ and:

    \[ \begin{array}[b]{rcl}
      \bivbinomial[1](\gamma)
      & = &
      \bivbinomial[1](\gamma)(0,0)\bigket{0,0} +
      \bivbinomial[1](\gamma)(0,1)\bigket{0,1}
      \\
      & & \quad + \,
      \bivbinomial[1](\gamma)(1,0)\bigket{1,0} +
      \bivbinomial[1](\gamma)(1,1)\bigket{1,1}
      \\
      & = &
      \multinomial[1](\gamma)\big(1\ket{0,0}\big)\bigket{0,0} +
      \multinomial[1](\gamma)\big(1\ket{0,1}\big)\bigket{0,1}
      \\
      & & \quad + \,
      \multinomial[1](\gamma)\big(1\ket{1,0}\big)\bigket{1,0} +
      \multinomial[1](\gamma)\big(1\ket{1,1}\big)\bigket{1,1}
      \\
      & = &
      \gamma(0,0)\bigket{0,0} + \gamma(1,0)\bigket{1,0} +
      \gamma(0,1)\bigket{0,1} + \gamma(1,1)\bigket{1,1}
      \hspace*{\arraycolsep}=\hspace*{\arraycolsep}
      \gamma.
    \end{array} \eqno{\QEDbox} \]
\end{myproof}

\end{document}